\documentclass[11pt,leqno]{amsart}
\usepackage{amssymb}
\usepackage{array}
\usepackage{tabularx}

\newcommand{\mesh}{\mathbin{\#}}

\newcommand{\cl}{\operatorname{cl}}

\newcommand{\adh}{\operatorname{adh}}

\newcommand{\lm}{\operatorname{lim}}

\def\A{{\mathcal A}}  \def\B{{\mathcal B}} \def\C{{\mathcal C}} \def\D{{\mathcal D}}
\def\E{{\mathcal E}}
\def\F{{\mathcal F}} \def\G{{\mathcal G}} \def\H{{\mathcal H}} 
\def\J{{\mathcal J}} \def\K{{\mathcal K}} \def\L{{\mathcal L}}
 \def\N{{\mathcal N}}
\def\O{{\mathcal O}} \def\P{{\mathcal P}}  \def\Q{{\mathcal Q}}
  \def\S{{\mathcal S}} 
\def\U{{\mathcal U}}

   \def\Fo{{\mathbb F}}

\def\then{\Longrightarrow\ }

\def\0{\varnothing}
\def\inc{\subset }
\def\iff{\ \Longleftrightarrow \ }

\def\to{\mathop{\rightarrow}\limits}

\newtheorem{prop}{Proposition}
\newtheorem{thm}[prop]{Theorem}
\newtheorem{lem}[prop]{Lemma}
\newtheorem{cor}[prop]{Corollary}

\theoremstyle{definition}

\newtheorem{ex}[prop]{Example}

\newtheorem{question}[prop]{Question}

\begin{document}
\title[Compatible relations]{Compatible relations on filters
and stability of local topological properties under supremum and product}
\author{Francis Jordan}
\author{Fr\'ed\'eric Mynard}
\address{Department of Mathematical Sciences, Georgia Southern University,
0203 Georgia Ave. Room 3008 Statesboro, GA 30460-8093}
\email{fmynard@GeorgiaSouthern.edu}
\email{fjordan@GeorgiaSouthern.edu}
\thanks{The authors want to thank S. Dolecki for comments on \cite{JM} 
which are at the origin of the approach in terms of relations developed in 
this paper.}
\date\today

\begin{abstract}
An abstract scheme using particular types of relations on filters
leads to general unifying results on stability under supremum and product
 of local topological properties. We present applications for Fr\'echetness,
strong Fr\'echetness, countable tightness and countable fan-tightness,
some of which recover or refine classical results, some of which are new.
The reader may find other applications as well.
\end{abstract}
\maketitle
\section{Introduction}
A large number of topological properties fail to  be stable under
(even finite) product or even by supremum of topologies. Among
such properties are a lot of fundamental local topological
features such as Fr\'echetness, strong Fr\'echetness, countable
tightness and countable fan-tightness, to cite a few (see the next
section for definitions). Consequently, 
the quest for conditions on the factor spaces to ensure that the
product space has the desired local property has attracted a lot
of attention, e.g., \cite {Ar.spectrum}, \cite{cost.simon},
\cite{dolnog}, \cite{vandowen}, \cite{gerlits}, \cite{gruenhage},
\cite{Kendrick}, \cite{quest}, \cite{Nogura87}, \cite{Nogura88},
\cite {Novak1}, \cite{Novak2}, \cite{Novak3}, \cite{nyikos},
\cite{olson}, \cite{simon}, \cite{tamano}, \cite{JM},
\cite{JMCRAS} for Fr\'echetness and strong Fr\'echetness alone;
\cite{Ar.spectrum}, \cite{bella.malykhin}, \cite{ar.bella.CMUC},
\cite{bella.vanmill}, \cite{CHV}, \cite{quest}, \cite{malyhin} for
tightness and fan-tightness.

In this paper, we propose a unified approach to this class of
problems, obtaining as byproducts of our theory both refinements
and generalizations of known results and entirely new theorems.
More specifically, we are interested in the following type of
problem: Let $(\mathcal{P})$ and $(\mathcal{Q})$ be local
topological properties.
\begin{question}\label{prob1}
Characterize topological spaces $X$ such that $X\times Y$ has
property $(\mathcal Q)$ for every space $Y$ with property
$(\mathcal P)$.
\end{question}
In most cases, investigations have been restricted to
$\mathcal{(P)=(Q)}$.

The first crucial observation is that if $\P$ and $\Q$ are local
topological properties, they are characterized by a property of
the neighborhood filters. Of course, the property is not stable
under product because the corresponding class of filters is not.
The second crucial observation is that, even if these classes of
filters behave badly with respect to the product operation, they
are almost always defined via other better behaved classes of
filters.  For instance,  a topological space is {\em Fr\'echet} if whenever
a point $x$ is in the closure of a subset $A$, there exists a
sequence (equivalently a countably based filter) on $A$ that
converges to $x$. 

If $A\subseteq X$, then the {\em principal filter of} $A$ is $A^\uparrow=\{B\subseteq
X: A\subseteq B\}$. Similarly, if $\A\inc 2^X$, we denote $\A^\uparrow$ the family of subsets of $X$ that contains an element of $\A$. However,
we will often identify subsets of $X$ with their
principal filters, that is,  $A^\uparrow$ is often simply denoted $A$. 
 Two collections of sets $\F$ and $\G$ {\em mesh},
in symbol $\F\mesh\G$, if $F\cap G\neq \varnothing$ for every $F\in
\F$ and every $G\in \G$. The supremum $\F\vee\G$ of two
collections of sets $\F$ and $\G$ exists when $\F\mesh\G$ and stands
for the collection of intersections $\{F\cap G\colon F\in\F\text{
and }G\in\G\}$.  However, when $\F$ and $\G$ are filters we will
understand $\F\vee\G$ to be the filter $(\F\vee\G)^\uparrow$ generated by the collection
of intersections.  

With these notations, the definition of
Fr\'echetness in terms of closure rephrases in terms of
neighborhood filter as follows: $\F$ is a Fr\'echet filter (a
filter of neighborhood in a Fr\'echet space) if whenever $A\mesh\F$,
there exists a countably based filter (equivalently, a sequence)
$\L$ which is finer than $A\vee\F$. This definition depends on the
class $\mathbb {F}_1$ of principal filters and on the class
$\mathbb {F}_\omega$ of countably based filters, which are both
productive.

We will take an approach based on relations between filters.  
For example, we consider the relation $\bigtriangleup$ on the
set $\mathbb {F}(X)$ of filters on $X$ defined by
$\F\bigtriangleup \H$ if
$$\H\mesh\F \then \exists \L \in \mathbb{F}_\omega: \L\geq \F\vee\H,$$
a filter is Fr\'echet if and only if it is in the
$\bigtriangleup$-relation with every principal filter. In other
words the class of Fr\'echet filters is $\mathbb
{F}_1^\bigtriangleup$ \cite{malyhin}, \cite{dol.active} where we denote by $\mathbb {J}^\star$ the
filters that are in the $\star$-relation with every filter of the
class $\mathbb {J}$.

Under mild conditions on $\mathbb {J}$ and on relations $\star$
and $\square$, we characterize filters whose supremum (in section
3)
 and whose product (section 4) with a filter of the class $\mathbb {J}$ or
 $\mathbb {J}^\star$ is a filter of $\mathbb {J}^\star$ or $\mathbb {J}^\square$. The quests
 for stability under supremum and under product turn out to be
 intimately related (section 4). These abstract results leads to a
 large collection of significant concrete corollaries, because
 most classical local topological properties, like  Fr\'echetness, can be characterized
 in terms of neighborhood filters of the type $\mathbb {J}^\star$ for
 classes $\mathbb {J}$ and relations $\star$ that fulfill the needed
 conditions.

In section 4, we show how solutions for the problem of stability
 under product are obtained as an abstractly defined subclass, called kernel, of the
 class of solutions for the problem of stability under supremum.
 From the technical viewpoint, the difficulties (that attracted attention to this
 type of problems for so long) lie in the internal
 characterization of kernels. In section 5, we characterize a variety
 of such kernels, obtaining as byproducts improvements of
 classical results as well as entirely new results.
 Finally, our abstract approach allows us to clarify (section 6)
 the relationships between all these properties, improving again upon
  the  known results.

We present applications related to three relations only, but the
theory is designed for further applications with other examples of
relations and of classes of filters.

\section{sup-compatible relations on $\mathbb {F}(X)$.}
Let $R$ denote a relation on a set $X$ ($R\subseteq X\times X$). As usual,
$R x=\{y\in X: (x,y)\in R\}$ and if $F\subseteq X$,
$R F = \bigcup_{x \in F} R x$.
The {\em polar of } a subset $F$ of $X$ (with respect to $R$) is
$$F^R = \bigcap_{x\in F} R x,$$
with the convention that $\varnothing^R=X$.
An immediate consequence of the definitions is that
\begin{lem} \label{lem:polar}
If $R$ is a symmetric relation on $X$ and $F\subseteq X$, then
\begin{enumerate}
\item $F\subseteq F^{RR}$;
\item $F^{RRR}=F^{R}$.
\end{enumerate}
Moreover, if $F\subseteq F^{R}$ then $F\subseteq F^{RR}\subseteq F^{R}$.
\end{lem}

The symbol $\neg$ denotes negation. For instance, $x(\neg R) y$ means that $(x,y)\notin R$.

The set of filters on a given set $X$ is denoted by $\mathbb {F}
(X)$. A symmetric relation $\star$ on  $\mathbb {F}(X)$ for which
$\F\star\G$ whenever $\F$ and $\G$ do not mesh and which verifies
\begin{equation}\tag{$\bullet$}
F\star (\G\vee\H) \then (\F\vee\H)\star \G,
\end{equation}
whenever $\G\mesh\F$ and $\F\mesh\H$ is called {\em sup-compatible} or
$\vee$-{\em compatible}.

We are going to consider several particular classes of filters. In
general $\mathbb {J}$ and $\mathbb {D}$ will denote generic
classes of filters. A filter of the class $\mathbb {J}$ is a
called $\mathbb {J}$-{\em filter}. The collection of $\mathbb
{J}$-filters on $X$ is denoted $\mathbb {J}(X)$. However, we will
often omit $X$ when the underlying set considered is clear.
 In particular, we will frequently consider the classes
$\mathbb {F}_1$ and $\mathbb {F}_{\omega}$ of principal and
countably based filters, respectively.

\begin{ex}  \label{ex:frechet}
The relation $\bigtriangleup$ on $\mathbb {F}(X)$ defined by $\F\bigtriangleup \H$ if
\begin{equation}\label{frechet}
\H\mesh\F \then \exists \L \in \mathbb{F}_\omega: \L\geq \F\vee\H,
\end{equation}
 is a $\vee$-compatible relation.  Then $\mathbb F_1^\bigtriangleup$ is the class of {\em Fr\'echet filters} \cite{malyhin},
\cite{dol.active}. As noticed in the introduction,
a topological space is Fr\'echet if and only if all its
neighborhood filters are Fr\'echet.

Analogously, $\mathbb{F}_\omega^\bigtriangleup$ is the class of
{\em strongly Fr\'echet filters} \cite{malyhin},
\cite{dol.active}, that is, of filters $\F$ satisfying (\ref{frechet}) for every
countably based filter $\H$.
Recall that a topological space $X$ is {\em strongly Fr\'echet} if for every
$x\in X$ and every decreasing countable collection $A_n\subseteq X$ such that
$x\in \cl(A_n)$ for every $n$, there exists a sequence $x_n \in A_n$ that converges to $x$.
It is easy to see that a topological space is
strongly Fr\'echet if and only if all its neighborhood filters are strongly Fr\'echet.
\end{ex}

\begin{ex}  \label{ex:tight}
Consider the relation $\diamondsuit$ on $\mathbb {F}(X)$ defined
by $\F{\diamondsuit}\H$ if
$$
\H\mesh\F \then 
\exists A \in \mathbb {F}_1, |A|\leq\omega: A\mesh \F\vee\H.
$$
This is a $\vee$-compatible relation. Then $\mathbb
{F}_1^\diamondsuit$ is the class of {\em countably tight filters}.
Recall that a topological space $X$ is countably tight
\cite{Ar.spectrum} if for every $x\in X$ and $A\subseteq X$ such
that $x\in \cl(A)$,
 there exists a countable subset
$B\subseteq A$ such that  $x\in \cl(B)$. It is easy to see that a space $X$ is countably tight  if and only if all its neighborhood filters are countably tight.
\end{ex}

\begin{ex} \label{ex:fantight}
A topological space is {\em countably fan-tight}
\cite{ar.bella.CMUC}
if for every countable family $(A_n)_\omega$ of
subsets such that $x\in \bigcap_{n\in\omega} \cl(A_n)$, 
there exists finite subsets
$B_n$ of $A_n$ such that $x\in \cl(\bigcup_{n\in\omega} B_n)$.

We call a filter $\F$ {\em countably fan-tight} if whenever $A_n\mesh\F$, 
there exists finite subsets $B_n$ of $A_n$ such
that $\bigcup_{n\in\omega} B_n\mesh\F$. 
It was observed in \cite[Remark 1]{ar.bella.CMUC} that the definition of 
fan-tightness is unchanged if we only consider {\em decreasing} countable 
collections $(A_n)_\omega$.
Clearly, a space is countably fan-tight if and only if all its neighborhood filters are 
countably fan-tight.
Consider the relation $\dagger$ on $\mathbb {F}(X)$ defined by
$\F \dagger \H$ if
$$(A_n)_\omega\mesh(\F\vee\H) \then 
\exists B_n\subseteq A_n, |B_n|<\omega,
(\bigcup_{n\in\omega} B_n)\mesh \F\vee\H.
$$
This is a $\vee$-compatible relation.
By \cite[Remark 1]{ar.bella.CMUC}, $\F \dagger  \H$ if for any {\em decreasing} countable filter base $(A_n)_\omega$ meshing with $\F\vee\H$, there exists finite sets $B_n\subseteq A_n$ such that 
$(\bigcup_{n\in\omega} B_n)\mesh \F\vee\H$.
\end{ex}

\begin{lem}\label{lem:fanT}
The following are equivalent:
\begin{enumerate}
\item $\F$ is countably fan tight;\label{cft1}
\item \label{cft2}
$\{E_k\}_{k\in\omega}\mesh\F$, then there exist finite sets
$B_k\subseteq E_k$ such that $\{\bigcup_{n\leq
k}B_k\}_{n\in\omega}\mesh\F$;
\item \label{cft3}
If $\{E_k\}_{k\in\omega}$ is a decreasing countable filter base and $\{E_k\}_{k\in\omega}\mesh\F$, 
then there exist finite sets $B_k\subseteq E_k$ such that $\{\bigcup_{n\leq
k}B_k\}_{n\in\omega}\mesh\F$;
\item $\F\in\mathbb {F}_\omega^\dagger$;\label{cft4}
\item $\F\in\mathbb {F}_1^\dagger$.\label{cft5}
\end{enumerate}
\end{lem}
\begin{proof}
(\ref{cft1}) $\then$ (\ref{cft2}). 
Fix $n\in\omega$.  Since $E_k\mesh\F$ for all $k\geq n$, there exist
finite sets $\{B^n_k\}_{k\geq n}$ such that $B_k^n\subseteq E_k$
and $(\bigcup_{k\geq n}B^n_k)\mesh\F$.
 For every $k\in\omega$ let $B_k=\bigcup_{k\geq n}B_k^n$.
 Notice $B_k\subseteq E_k$ and $B_k$ is finite for every $k\in\omega$.
Let $n\in\omega$ and $F\in\F$.
There is a $k\geq n$ such that $F\cap B^n_k\neq\emptyset$.
  Since $k\geq n$, we have $B^n_k\subseteq B_k$.
 Thus, $k\geq n$ and $B_k\cap F\neq\emptyset$.
Therefore, $(\bigcup_{n\leq k}B_k)\mesh\F$ for every $n\in\omega$.

(\ref{cft2}) $\then$ (\ref{cft3}) and (\ref{cft4}) $\then$ (\ref{cft5}) are straightforward.
(\ref{cft5}) $\then$ (\ref{cft1}) was observed in Example \ref{ex:fantight}. 

(\ref{cft3}) $\then$ (\ref{cft4}) follows from \cite[Remark 1]{ar.bella.CMUC}. Indeed, if
 $\F$ is as in (\ref{cft3}), $\H$ and $\A$  are countably based filters with decreasing filter bases $(H_n)_{n\in\omega}$ and $(A_n)_{n\in\omega}$ respectively, such that
$(A_n)_n\mesh(\H\vee\F)$, then $(A_n\cap H_n)_n\mesh\F$. Therefore, there exists finite sets  $B_n\inc A_n\cap H_n$ such that  
$\{\bigcup_{n\leq k}B_k\}_{n\in\omega}\mesh\F$. Clearly, 
$\{\bigcup_{n\leq k}B_k\}_{n\in\omega}\mesh(\F\vee\H)$ and $\F\in\mathbb {F}_{\omega}^\dagger$.
\end{proof}

\section{Stability of local properties under supremum}
We call a class $\mathbb {J}$ of filters  $(\mathbb {D},\mathbb
{M})$-{\em steady} if $\F\vee\H\in \mathbb {M}(X)$ whenever
$\F\in\mathbb {J}(X)$ and $\H\in\mathbb {D}(X)$. If $\mathbb {J}$
is $(\mathbb {D},\mathbb {J})$-steady we say that $\mathbb {J}$ is
$\mathbb {D}$-{\em steady}. If $\mathbb {J}$ is $\mathbb
{J}$-steady, we simply say that $\mathbb {J}$ is {\em steady}. By
Lemma \ref{lem:polar}, we only generate the two classes $\mathbb
{J}^\star$ and $\mathbb {J}^{\star\star}$ from a given class
$\mathbb {J}$ by taking the polars with respect to a symmetric
relation $\star$ on filters. We now investigate stability
relationships between these classes under supremum.  To begin
notice that an immediate consequence of the definitions is that a
class of filters $\mathbb {J}$ is $(\mathbb {D},\mathbb {M})$-{\em
steady} if and only if $\mathbb {D}$ is $(\mathbb {J},\mathbb
{M})$-{\em steady}.
\begin{prop}    \label{prop:stability}
If $\mathbb {J}$ is a $\mathbb D$-steady class of filters, and if
$\star$ and $\square$ are $\vee$-compatible relations,
 then $\mathbb {J}^\star$ and $\mathbb {J}^{\star\square}$ are both $\mathbb {D}$-steady.
\end{prop}
\begin{proof}
Let $\F\in \mathbb {J}^\star(X)$ and $\H\in \mathbb {D}(X)$ such
that $\F\mesh\H$. We want to show that $\F\vee\H\in\mathbb
{J}^\star$. Let $\L$ be a $\mathbb {J}$-filter such that
$\L\mesh\F\vee\H$. As $\mathbb {J}$ is $\mathbb D$-steady,
$\L\vee\H\in \mathbb {J}$. Thus, $\F\star (\L\vee\H)$ because
$\F\in\mathbb {J}^\star$. By $(\bullet)$, $(\F\vee\H)\star \L$.
Hence, $\F\vee\H\in\mathbb {J}^\star$.

Let $\F\in \mathbb {J}^{\star\square}(X)$ and $\H\in \mathbb
{D}(X)$ such that $\F\mesh\H$. We want to show that
$\F\vee\H\in\mathbb {J}^{\star\square}$. Let $\L$ be a $\mathbb
{J}^\star$-filter such that $\L\mesh\F\vee\H$. As we have proved that
$\mathbb {J}^\star$ is $\mathbb {D}$-steady, $\L\vee\H\in \mathbb
{J}^\star$. Thus, $\F\square (\L\vee\H)$ because $\F\in\mathbb
{J}^{\star\square}$. By ($\bullet$), $(\F\vee\H)\square \L$.
Hence, $\F\vee\H\in\mathbb {J}^{\star\square}$.
\end{proof}

\begin{cor} \label{cor:comp}
If $\mathbb {J}$ is a steady class of filters and $\star$ and
$\square$ are $\vee$-compatible relations, then $\mathbb
{J}^\star$ is $(\mathbb {J}^{\star\square},\mathbb
{J}^\square)$-steady. In particular, $\mathbb {J}^\star$ is
$\mathbb {J}^{\star\star}$-steady.
\end{cor}
\begin{proof}
Let $\F\in \mathbb {J}^{\star}(X)$ and $\H\in \mathbb
{J}^{\star\square}(X)$ such that $\F\mesh\H$. We want to show that
$\F\vee\H\in\mathbb {J}^{\square}$. Let $\L$ be a $\mathbb
{J}$-filter such that $\L\mesh\F\vee\H$. As $\mathbb
{J}^{\star\square}$ is $\mathbb {J}$-steady, $\L\vee\H\in \mathbb
{J}^{\star\square}$, so that $(\L\vee\H)\square\F$. By
($\bullet$), $(\F\vee\H)\square\L$.
\end{proof}

\begin{cor} \label{cor:comp1}
If $\mathbb {J}$ is a steady class of filters and $\star$,
$\square$, and $\bigtriangledown$ are $\vee$-compatible relations,
then $\mathbb {J}^{\square\bigtriangledown}$ is $(\mathbb
{J}^{\star\square},\mathbb {J}^{\star\bigtriangledown})$-steady.
In particular, $\mathbb {J}^{\star\star}$ is steady.
\end{cor}
\begin{proof}
Let $\F\in \mathbb {J}^{\square\bigtriangledown}(X)$ and $\H\in
\mathbb {J}^{\star\square}(X)$ such that $\F\mesh\H$. We want to show
that $\F\vee\H\in\mathbb {J}^{\star\bigtriangledown}$. Let $\L$ be
a $\mathbb {J}^{\star}$-filter such that $\L\mesh\F\vee\H$. As
$\mathbb {J}^{\star}$ is $(\mathbb
{J}^{\star\square},\mathbb{J}^{\square})$-steady, $\L\vee\H\in
\mathbb {J}^{\square}$, so that $(\L\vee\H)\bigtriangledown\F$. By
($\bullet$), $(\F\vee\H)\bigtriangledown\L$.
\end{proof}


\begin{thm} \label{th:main}
Let $\star$ and $\square$ be two $\vee$-compatible relations on
$\mathbb {F}(X)$ and let $\mathbb {J}$ be a steady class of
filters containing $\mathbb {F}_1$.
\begin{enumerate}
\item\label{it01}
\begin{eqnarray*}
\F\in \mathbb {J}^\star &\iff& \forall \G\in\mathbb {J}, \G\mesh\F\then \F\vee \G\in \mathbb {F}_1^\star\\
&\iff& \forall \G\in \mathbb {J}, \G\mesh\F\then \F\vee\G\in \mathbb {J}^\star\\
&\iff& \forall \G\in \mathbb {J}^{\star\star}, \G\mesh\F\then
\F\vee\G\in \mathbb {J}^\star.
\end{eqnarray*}
\item \label{it02}
\begin{eqnarray*}
\F\in \mathbb {J}^{\star\square} &\iff& \forall \G\in \mathbb
{J}^\star,
\G\mesh\F\then \F\vee\G\in \mathbb {J}^{\star\star\square}\\
&\iff& \forall \G\in \mathbb {J}^\star,
\G\mesh\F\then \F\vee\G\in \mathbb {F}_1^{\star\star\square}\\
&\iff& \forall \G\in \mathbb {J}^\star,
\G\mesh\F\then \F\vee\G\in \mathbb {J}^\square\\
&\iff& \forall \G\in \mathbb {J}^\star, \G\mesh\F\then \F\vee\G\in
\mathbb {F}_1^\square.
\end{eqnarray*}
\end{enumerate}

\end{thm}
\begin{proof}
We show (\ref{it01}).  If $\F\in\mathbb {J}^\star$ then
$$\forall \G\in \mathbb {J}^{\star\star}, \G\mesh\F\then \F\vee\G\in \mathbb {J}^\star$$
because $\mathbb {J}^\star$ is $\mathbb {J}^{\star\star}$-steady,
by Corollary \ref{cor:comp}. The two other direct implications
follow from $\mathbb {F}_1\subseteq\mathbb {J}\subseteq\mathbb
{J}^{\star\star}$ (the last inclusion
 comes from Lemma \ref{lem:polar}).

Conversely, if $\F\notin \mathbb {J}^\star(X)$, then there exists
$\G\in \mathbb {J}(X)$ such that $\F(\neg\star)\G$, or in other
words, $\F(\neg \star)(\G\vee X)$. Therefore $(\F\vee \G) (\neg
\star) X$, so that $(\F\vee\G) \notin \mathbb {F}_1^\star$.
We show (\ref{it02}).  Let  $\F\in\mathbb {J}^{\star\square}$ and
$\G\in\mathbb {J}^{\star}$.  Let $\H\in\mathbb {J}^{\star\star}$
and $\H\mesh(\F\vee\G)$.  Since $\H\mesh\G$, we have, by
Corollary~\ref{cor:comp}, $\H\vee\G\in\mathbb {J}^{\star}$.  So,
$\F\square(\H\vee\G)$.  By ($\bullet$), $(\F\vee\G)\square\H$. So,
$\F\vee\G\in\mathbb {J}^{\star\star\square}$.  By the containments
$\mathbb {J}^{\star\star\square}\subseteq\mathbb
{J}^{\square}\subseteq\mathbb {F}_1^{\square}$ and $\mathbb
{J}^{\star\star\square}\subseteq\mathbb
{F}_1^{\star\star\square}\subseteq\mathbb {F}_1^{\square}$,
$\F\in\mathbb {J}^{\star\square}$ implies any of the first three
statements on the right and each of the first three statements on
the right implies the last statement on the right.

Conversely, assume that $\F\notin \mathbb {J}^{\star\square}$.
Then, there exists $\H\in\mathbb {J}^\star$ such that $\H\mesh\F$ but
$\H(\neg\square) \F$. By ($\bullet$), $\H\vee X (\neg \square)
\F$ and $\F\vee\H(\neg \square) X$. Thus, $\F\vee\H\notin
\mathbb {F}_1^\square$.
\end{proof}

A topological space is called $\mathbb {J}$-{\em based} if all its
neighborhood filters are $\mathbb {J}$-filters. For instance, the
$\mathbb{J}$-based spaces are respectively the {\em finitely
generated} \cite{lowen.sonck}\cite{myn.cont}, first-countable,
Fr\'echet, strongly Fr\'echet, countably tight, countably
fan-tight spaces, when $\mathbb{J}$ is  the class $\mathbb{F}_1$,
$\mathbb{F}_\omega$, $\mathbb{F}_1^\bigtriangleup$,
$\mathbb{F}_\omega^\bigtriangleup$, $\mathbb{F}_1^\diamondsuit$
and $\mathbb{F}_1^\dagger$ respectively.
 We say that a class
$\mathbb {J}$ of filters is called {\em pointable} if $\F\wedge\{x\}
\in \mathbb {J}(X)$ whenever $\F\in\mathbb {J}(X)$, for every set
$X$ and every $x\in X$.
\begin{cor} \label{cor:sup}
Let $\star$, $\bigtriangledown$, and $\square$ be
$\vee$-compatible relations on filters. Let $\mathbb {J}$ be a
pointable and steady class containing $\mathbb {F}_1$.
\begin{enumerate}
\item Each of the following statements are equivalent:
\begin{enumerate}
\item \label{it11}
$(X,\tau)$ is $\mathbb {J}^\star$-based,
\item \label{it14}
$(X,\tau\vee\xi)$ is $\mathbb {J}^\star$-based for every $\mathbb
{J}^{\star\star}$-based topology $\xi$ on $X$,
\item \label{it13}
$(X,\tau\vee\xi)$ is $\mathbb {J}^\star$-based for every $\mathbb
{J}$-based topology $\xi$ on $X$,
\item \label{it12}
$(X,\tau\vee\xi)$ is $\mathbb {F}_1^\star$-based for every
$\mathbb {J}$-based topology $\xi$ on $X$.
\end{enumerate}
\item
If $\mathbb {J}^\star$ is pointable then the following are
equivalent:
\begin{enumerate}
\item \label{it21}
$(X,\tau)$ is $\mathbb {J}^{\star\square}$-based,
\item\label{it20}
$(X,\tau\vee\xi)$ is $\mathbb {J}^{\star\star\square}$-based for
every $\mathbb {J}^\star$-based topology $\xi$ on $X$,
\item\label{it24}
$(X,\tau\vee\xi)$ is $\mathbb {F}_1^{\star\star\square}$-based for
every $\mathbb {J}^\star$-based topology $\xi$ on $X$,
\item \label{it23}
$(X,\tau\vee\xi)$ is $\mathbb {J}^\square$-based for every
$\mathbb {J}^\star$-based topology $\xi$ on $X$,
\item \label{it22}
$(X,\tau\vee\xi)$ is $\mathbb {F}_1^\square$-based for every
$\mathbb {J}^\star$-based topology $\xi$ on $X$.
\end{enumerate}
\item
If $(X,\tau)$ is $\mathbb {J}^{\square\bigtriangledown}$-based, then
$(X,\tau\vee\xi)$ is $\mathbb {J}^{\star\bigtriangledown}$-based
for every $\mathbb {J}^{\star\square}$-based topology $\xi$ on
$X$.
\end{enumerate}
\end{cor}
\begin{proof}
(\ref{it11}) $\then$ (\ref{it14}) follows immediately from Theorem
\ref{th:main}. (\ref{it14}) $\then$ (\ref{it13}) $\then$
(\ref{it12}) follows from $\mathbb {J}\subseteq\mathbb
{J}^{\star\star}$ and $\mathbb {J}^\star\subseteq\mathbb
{F}_1^\star$ (because $\mathbb {F}_1\subseteq \mathbb {J}$).

(\ref{it12}) $\then$ (\ref{it11}). By way of contradiction,
assuming that $(X,\tau)$ is not $\mathbb {J}^\star$-based, there
exists $x\in X$ such that $\N_\tau(x)\notin \mathbb {J}^\star$. In
view of the first part of Theorem \ref{th:main}, there exists
$\G\in \mathbb {J}$ such that $\G\mesh\N_\tau(x)$ and
$\N_\tau(x)\vee\G \notin \mathbb {F}_1^\star$. Let $\xi$ be the
topology on $X$ with all points but $x$ isolated with
$\N_\xi(x)=\G\wedge(x)$. The space $(X,\xi)$ is $\mathbb
{J}$-based because $\mathbb {J}$ is pointable. So, by
(\ref{it12}),
$\N_\tau(x)\vee(\G\wedge(x))=\N_{\tau\vee\xi}(x)\in\mathbb
{F}_1^\star$. We consider two cases.

If $x\notin\bigcap\G$, then $(\N_{\tau}(x)\vee(\G\wedge(x)))\mesh(X\setminus\{x\})$ and
$\G=(\G\wedge(x))\vee(X\setminus\{x\})$.  So, $(\N_{\tau}(x)\vee(\G\wedge(x)))\star(X\setminus\{x\})$.
 Since $\star$ is $\vee$-compatible, $\N_{\tau}(x)\star(((\G\wedge(x))\vee(X\setminus\{x\}))$.
 Thus, $\N_{\tau}(x)\star\G$, a contradiction.

If $x\in\bigcap\G$, then $\G\wedge(x)=\G$.  In this case we have
$\N_{\tau}(x)\vee\G\in\mathbb {F}_1^\star$ and $\N_{\tau}(x)\mesh\G$.
Since $(\N_{\tau}(x)\vee\G)\mesh X$, $(\N_{\tau}(x)\vee\G)\star X$. By
the $\vee$-compatibility of $\star$, $\N_{\tau}\star(\G\vee X)$.
Thus, $\N_{\tau}(x)\star\G$, a contradiction.

Therefore, $(X,\tau)$ is $\mathbb {J}^\star$-based.

The second part is proved in a similar way.

The third part follows from Corollary~\ref{cor:comp1}.

%
\end{proof}

A class $\mathbb {J}$ of filters is called $\mathbb {F}_1$-{\em
composable} if the image of a $\mathbb {J}$-filter under a
relation is a $\mathbb {J}$-filter.
\begin{lem}     \label{lem:phi1}
\begin{enumerate}
\item $\mathbb {F}_1$, $\mathbb {F}_\omega$, $\mathbb {F}_1^{\diamondsuit}$,
$\mathbb {F}_1^\dagger$, $\mathbb {F}_1^\bigtriangleup$, $\mathbb
{F}_\omega^\bigtriangleup$ are all $\mathbb {F}_1$-composable.
\item an $\mathbb {F}_1$-composable class is pointable.
\end{enumerate}
\end{lem}
\begin{proof}
For brevity we only show that $\mathbb {F}_1^{\dagger}$ is
$\mathbb {F}_1$-composable, the other cases are similar and more
straightforward.  Let $\F\in\mathbb {F}_1^{\dagger}(X)$, $Y$ be a
set and $R\subseteq X\times Y$.  Suppose $B\subseteq Y$ and
$(A_n)_{\omega}\mesh(B\vee R\F)$.   Since
$R^{-}(A_n)_{\omega}\mesh(\F\vee R^{-}B)$, there exist finite sets
$K_n\subseteq R^{-}A_n$ such that $(\bigcup_{\omega}K_n)\mesh(\F\vee
R^{-}B)$.  For each $n$ pick a finite $J_n\subseteq A_n\cap B$
such that $K_n\subseteq R^{-}J_n$.  Let $F\in \F$.  There is an
$n$ such that $K_n\cap F\neq\emptyset$.  Since $K_n\subseteq
RJ_n$, $RJ_n\mesh F$.  So, $J_n\cap RF\neq\emptyset$.  Thus,
$(\bigcup_{\omega}J_n)\mesh(B\vee R\F)$ and $R\F\dagger B$.

For the second part consider the relation $R=\{(x,x)\colon x\in X\}\cup(X\times\{p\})$. Then $\F\wedge(p)=R\F$.
\end{proof}

We call  a filter $\F$ {\em almost principal} if there exists
$F_0\in\F$ such that $|F_0\setminus F|<\omega$ for every $F\in\F$. Principal filters and cofinite filters (of an infinite set) are almost principal.
We use the same name for spaces based in such filters.  Such
spaces include sequences and one point compactifications of
discrete sets. It is easily verified that every almost principal
filter is Fr\'echet.

If $\H\in \mathbb F(X)$, we denote by $\H^\bullet$ the {\em principal part} $\bigcap\H$ of $\H$ and by $\H^\circ$ the {\em free part} $\H\vee(\H^\bullet)^c$. One or the other may be the degenerate filter $\{\varnothing\}^\uparrow$. With the convention, that $\G\wedge\{\varnothing\}^\uparrow=\G$ for any filter $\G$, we have
\begin{equation*}
\H=\H^\circ\wedge\H^\bullet.
\end{equation*}


\begin{thm} \label{th:almostprinc}
$\mathbb {F}_1^{\bigtriangleup\bigtriangleup}$ is exactly the
class of almost principal filters.
\end{thm}
\begin{proof}
Let $\F$ be an almost principal filter and let $\H\mesh\F$ be a Fr\'echet filter.
There exists $F_0\in\F$ such that $|F_0\setminus F|<\omega$ for every $F\in\F$.
If $\H^\bullet(\neg\mesh)\F$, then $\H^{\circ}\mesh\F$.  In particular, $\H^{\circ}\mesh F_0$.  
So, there exists a free sequence finer than
$\H\vee F_0$, by Fr\'echetness
of $\H$. This sequence is also finer than $\F$, hence finer than $\F\vee\H$.

If $\H^{\bullet}\mesh\F$, then since $\F$ is Fr\'echet, there is a sequence finer than 
$\H^{\bullet}\vee\F\geq\H\vee\F$.  Thus, $\F\in \mathbb
{F}_1^{\bigtriangleup\bigtriangleup}$.

Conversely,
if $\F$ is not almost principal, then for all $F\in \F$ there exists $H_F\in \F$ such
that $|F\setminus H_F|\geq \omega$. Therefore, there exists a free sequence
$(x_n^F)_n$ on $F\setminus H_F$.
The filter $\bigwedge_{F\in \F} (x_n^F)_n$ is a Fr\'echet filter
 meshing with $\F$.  If $(y_n)_n$ is finer than $\F$, then for every
$F\in \F$, there exists $k_F$ such that $\{y_n\colon n\geq
k_F\}\subseteq H_F$. Therefore, there exists $n_F$ such that
$\{x_n^F\colon n\geq n_F\}\cap \{y_n\colon n\in
\omega\}=\emptyset$. The set $\bigcup_{F\in \F} \{x_n^F: n\geq n_F\}$ is 
an element of $\bigwedge_{F\in \F} (x_n^F)_n$
disjoint from $\{y_n:n\in \omega\}$. Thus, $\F\notin \mathbb
{F}_1^{\bigtriangleup\bigtriangleup}$.
\end{proof}

When $\star=\square=\bigtriangleup$ and $\mathbb{J=F}_1$,
Corollary \ref{cor:sup} rephrases as
\begin{cor} \label{cor:supFrechet}
\begin{enumerate}
\item
The following are equivalent:
\begin{enumerate}
\item $(X,\tau)$ is Fr\'echet;
\item $(X,\tau\vee\xi)$ is Fr\'echet for every finitely generated topology
$\xi$ on $X$;
\item $(X,\tau\vee\xi)$ is Fr\'echet for every
almost principal topology $\xi$ on $X$;
\end{enumerate}
\item
The following are equivalent:
\begin{enumerate}
\item $(X,\tau)$ is almost principal;
\item $(X,\tau\vee\xi)$ is Fr\'echet for every
Fr\'echet topology $\xi$ on $X$.
\end{enumerate}
\end{enumerate}
\end{cor}
No general condition was known to ensure that the supremum of two Fr\'echet topology
is Fr\'echet, as noticed for instance in \cite{CV}, \cite{CHV}.


Recall that $\mathbb {F}_\omega^{\bigtriangleup}$ is the class of
strongly Fr\'echet filters (or of neighborhood filters of strongly
Fr\'echet spaces).
 We call {\em productively Fr\'echet} \cite{JM} the filters from
the class $\mathbb {F}_\omega^{\bigtriangleup\bigtriangleup}$ and
we use  the same name for spaces based in such filters. When
$\star=\square=\bigtriangleup$ and $\mathbb{J=F}_\omega$,
Corollary \ref{cor:sup} rephrases as
\begin{cor} \label{cor:supstrFrechet}
\begin{enumerate}
\item
The following are equivalent:
\begin{enumerate}
\item $(X,\tau)$ is strongly Fr\'echet;
\item $(X,\tau\vee\xi)$ is Fr\'echet for every first-countable topology
$\xi$ on $X$;
\item $(X,\tau\vee\xi)$ is strongly Fr\'echet for every first-countable topology
$\xi$ on $X$;
\item $(X,\tau\vee\xi)$ is strongly Fr\'echet for every
productively Fr\'echet topology $\xi$ on $X$;
\end{enumerate}
\item
The following are equivalent:
\begin{enumerate}
\item $(X,\tau)$ is productively Fr\'echet;
\item $(X,\tau\vee\xi)$ is Fr\'echet for every
strongly Fr\'echet topology $\xi$ on $X$;
\item $(X,\tau\vee\xi)$ is strongly Fr\'echet for every
strongly Fr\'echet topology $\xi$ on $X$.
\end{enumerate}
\end{enumerate}
\end{cor}

Consider $\mathbb {J}=\mathbb {F}_1$ and $\star=\diamondsuit$. We
call {\em steadily countably tight} filters of $\mathbb
{F}_1^{\diamondsuit\diamondsuit}$ and we use the same name for
spaces based in such filters.  Notice that $\mathbb
{F}_1^{\diamondsuit}=\mathbb {F}_{\omega}^{\diamondsuit}$.
\begin{cor} \label{cor:suptight}
\begin{enumerate}
\item
The following are equivalent:
\begin{enumerate}
\item $(X,\tau)$ is countably tight;
\item $(X,\tau\vee\xi)$ is countably tight for every finitely
generated topology $\xi$ on $X$;
\item $(X,\tau\vee\xi)$ is countably tight for every countably based topology
$\xi$ on $X$;
\item $(X,\tau\vee\xi)$ is countably tight for every
steadily countably tight topology $\xi$ on $X$;
\end{enumerate}
\item
The following are equivalent:
\begin{enumerate}
\item $(X,\tau)$ is steadily countably tight;
\item $(X,\tau\vee\xi)$ is countably tight for every
countably tight topology $\xi$ on $X$.
\end{enumerate}
\end{enumerate}
\end{cor}

Consider $\mathbb {J}=\mathbb {F}_1$ and $\star=\dagger$. We call
{\em steadily countably fan-tight} filters of $\mathbb
{F}_1^{\dagger\dagger}$ and we use the same name for spaces based
in such filters.  Recall from Lemma~\ref{lem:fanT} that $\mathbb
{F}_1^\dagger=\mathbb {F}_{\omega}^\dagger$.

\begin{cor} \label{cor:supfantight}
\begin{enumerate}
\item
The following are equivalent:
\begin{enumerate}
\item $(X,\tau)$ is countably fan-tight;
\item $(X,\tau\vee\xi)$ is countably fan-tight for every finitely
generated topology $\xi$ on $X$;
\item $(X,\tau\vee\xi)$ is countably fan-tight for every countably based topology $\xi$ on $X$;
\item $(X,\tau\vee\xi)$ is countably fan-tight for every
steadily countably fan-tight topology $\xi$ on $X$;
\end{enumerate}
\item
The following are equivalent:
\begin{enumerate}
\item $(X,\tau)$ is steadily countably fan-tight;
\item $(X,\tau\vee\xi)$ is countably fan-tight for every
countably fan-tight topology $\xi$ on $X$.
\end{enumerate}
\end{enumerate}
\end{cor}
To our knowledge, no general conditions ensuring that the supremum of two countably (fan) tight topologies is countably (fan) tight (as in Corollaries \ref{cor:suptight} and \ref{cor:supfantight}) was known.

As a sample example of what one may get from Corollary~\ref{cor:sup} by
mixing relations we let $\mathbb {J}=\mathbb {F}_{\omega}$,
$\star=\bigtriangleup$,  $\square=\dagger$,
$\bigtriangledown=\diamondsuit$.

\begin{cor} \label{cor:strfrechettight}
\begin{enumerate}
\item
The following are equivalent:
\begin{enumerate}
\item $(X,\tau)$ is based in $\mathbb {F}_{\omega}^{\bigtriangleup\dagger}$;
\item $(X,\tau\vee\xi)$ is based in $\mathbb {F}_{\omega}^{\bigtriangleup\bigtriangleup\dagger}$
 for every strongly Fr\'echet topology $\xi$ on $X$;
\item $(X,\tau\vee\xi)$ is based in $\mathbb {F}_{1}^{\bigtriangleup\bigtriangleup\dagger}$
for every strongly Fr\'echet topology $\xi$ on $X$;
\item $(X,\tau\vee\xi)$ is countably fan tight for every strongly Fr\'echet topology
$\xi$ on $X$;
\end{enumerate}
\item If $(X,\tau)$ is based in $\mathbb {F}_{\omega}^{\bigtriangleup\dagger}$, then $(X,\tau\vee\xi)$
is based in $\mathbb {F}_{\omega}^{\bigtriangleup\diamondsuit}$
for every $\mathbb {F}_1^{\dagger\diamondsuit}$-based topology
$\xi$ on $X$.
\end{enumerate}
%

\end{cor}
\section{From steady to composable}

If $\F$ is a filter on a set $X$, $\G$ is a filter on a set $Y$ and $\H$ is a filter on
$X\times Y$ such that $\F\times Y \mesh\H$ and $X\times\G\mesh\H$, we denote by $\H\F$ the filter on $Y$ generated by the
sets $$HF=\{y: \exists x \in F, (x,y) \in H\},$$
for $H \in \H$ and $F \in \F$ and by $\H^-\G$ the filter on $X$ generated by
the sets
$$H^-G=\{x: \exists y \in G, (x,y) \in H\},$$
for $H\in \H$ and $G\in \G$. Notice that
$$\H\mesh(\F\times \G)\iff \H\F\mesh\G \iff \F\mesh\H^-\G.$$

A class $\mathbb {J}$ of filters is $(\mathbb {D},\mathbb
{M})$-{\em composable} if for every $X$ and $Y$ and every $\mathbb
{J}$-filter $\F$ on $X$ and every $\mathbb {D}$ filter $\H$ on
$X\times Y$, the filter $\H\F$ is an $\mathbb {M}$-filter on $Y$. A
class $\mathbb {J}$ is called $\mathbb {D}$-{\em composable} (\cite{DM}, \cite {myn.cont}) if it
is $(\mathbb {D},\mathbb {J})$-composable. We say that a class of
filters $\mathbb {J}$ is {\em projectable} provided that for every
$X$ and $Y$ and every $\mathbb {J}$-filter $\F$ on $X\times Y$ we
have $\pi_{Y}(\F)\in\mathbb {J}$.  Obviously, every $\mathbb
{F}_1$-composable class of filters is projectable.

Given a class of filters $\mathbb {J}$, we define the {\em kernel}
of $\mathbb {J}$ ($\ker(\mathbb {J})$) to be the class of
$(\mathbb {F}_1,\mathbb {J})$-composable filters.  Notice that
$\ker(\mathbb {J})$ is the largest $\mathbb {F}_1$-composable
subclass of $\mathbb {J}$.

\begin{lem}\label{lem:compstable}
If $\mathbb {D}$ is an $\mathbb {F}_1$-composable class of filters,
then $\mathbb {D}$ is $\mathbb {F}_1$-steady, and $\D\times A
\in\mathbb{D}$ for every $\D\in \mathbb{D}$ and every principal
filter $A$.
\end{lem}
\begin{proof}
Let $A\subseteq X$ and $\F\in\mathbb {D}(X)$ be such that $\F\mesh A$.
Since $\Delta_A=\{(x,x)\colon x\in A\}\in\mathbb {F}_1(X\times X)$
and $\mathbb {D}$ is $\mathbb {F}_1$-composable,
$A\vee\F=\Delta_A\F\in\mathbb {D}$.

Now, if $\D\in \mathbb{D}(X)$ and $A\inc Y$, then $\D\times
A=\pi_X^-\D\vee\pi_Y^-A \in \mathbb{D}$ because we already have
shown that $\mathbb{D}$ is $\mathbb{F}_1$-steady.
\end{proof}
\begin{thm}\label{th:eq}
\begin{enumerate}
\item \label{prodtosteady}
Let $\mathbb{M}$ be an $\mathbb{F}_1$-steady and projectable
class. If $\J\times \D\in \mathbb M$ for every $\J\in\mathbb{J}$
and every $\D\in\mathbb{D}$, then $\mathbb J$ is
$(\mathbb{D},\mathbb{M})$-steady.
\item \label{steadytocomposable}
Let $\mathbb{M}$ be projectable and $\mathbb J$ be
$\mathbb{F}_1$-composable. If $\mathbb J$ is
$(\mathbb{D},\mathbb{M})$-steady then $\mathbb J$ is
$(\mathbb{D},\mathbb{M})$-composable.
\item \label{composabletoproduct}
Let $\mathbb{D}$ be an $\mathbb{F}_1$-composable class. If
$\mathbb{J}$ is $(\mathbb{D},\mathbb{M})$-composable, then
$\J\times \D\in \mathbb M$ for every $\J\in\mathbb{J}$ and every
$\D\in\mathbb{D}$.
\item \label{eqker}
Let $\mathbb{M}$ be an $\mathbb{F}_1$-steady and projectable class.
If either $\mathbb{D}$ or $\mathbb{J}$ is
$\mathbb{F}_1$-composable then $\J\times \D\in \mathbb M$ for
every $\J\in\mathbb{J}$ and every $\D\in\mathbb{D}$ if and only if
$\J\times \D\in \mathbb \ker(\mathbb M)$ for every $\J\in\mathbb{J}$ and
every $\D\in\mathbb{D}$.
\end{enumerate}
\end{thm}
\begin{proof}
Proof of (\ref{prodtosteady}). Let $\D\in\mathbb D(X)$ and
$\J\in\mathbb J(X)$ such that $\D\mesh\J$. Then $\D\times
\J\in\mathbb M(X\times X)$. Let $\Delta$ be the diagonal of
$X\times X$. As $\mathbb M$ is $\mathbb F_1$-steady, $(\D\times
\J)\vee \Delta \in\mathbb M(X\times X)$. As $\mathbb M$ is
projectable, the $X$-projection of $(\D\times \J)\vee \Delta$ is an
$\mathbb M$-filter. We conclude with the observation that
$\pi_X((\D\times \J)\vee \Delta)=\D\vee\J$.

Proof of (\ref{steadytocomposable}). Let $\J\in\mathbb J(X)$ and let
$\H \in \mathbb D(X\times Y)$. The filter $\J\times Y$ is a
$\mathbb J$-filter  because $\mathbb J$ is $\mathbb
F_1$-composable. Therefore $(\J\times Y)\vee\H$ is an $\mathbb
M$-filter. As $\mathbb M$ is projectable, $\pi_Y((\J\times
Y)\vee\H)\in \mathbb M(Y)$. We conclude with the observation that
$\H(\J)=\pi_Y((\J\times Y)\vee\H)$.

Proof of (\ref{composabletoproduct}). Let $\J\in\mathbb J(X)$ and
$\D\in \mathbb D(Y)$. Let $\Delta_X$ be the diagonal of $X\times
X$, and consider the filter $\Delta_X\times \D$ as a filter of
relations from $X$ to $X\times Y$, that is, a filter on $X\times
(X\times Y)$. This is a $\mathbb D$-filter by Lemma
\ref{lem:compstable}, because $\mathbb D$ is $\mathbb
F_1$-composable. Moreover, $\J\times \D=(\Delta_X\times \D)(\J)$,
so that $\J\times \D\in\mathbb M(X\times Y)$.

Proof of \ref{eqker}. Let $\J\in\mathbb J(X)$, $\D\in\mathbb D(Y)$
and $A\inc X\times Y\times Z$. We want to show that $A(\J\times
\D)\in \mathbb M(Z)$. If either $\mathbb D$ or $\mathbb J$ is
$\mathbb F_1$-composable, then either $\J\times Z\in \mathbb J
(X\times Z)$ or $\D\times Z\in \mathbb D(Y\times Z)$. In any case,
$\J\times\D\times Z \in \mathbb M(X\times Y\times Z)$. Moreover,
$\pi_Z(A\vee(\J\times\D\times Z)) \in \mathbb M(Z)$ because
$\mathbb M$ is both  $\mathbb{F}_1$-steady and projectable. We
conclude with the observation that
$A(\J\times\D)=\pi_Z(A\vee(\J\times\D\times Z))$.
\end{proof}
The following diagram summarizes these relationships between
stability under product, composability and steadiness. Notice that
$\ker(\mathbb{M})$ is $\mathbb{F}_1$-composable, hence
$\mathbb{F}_1$-steady and projectable.

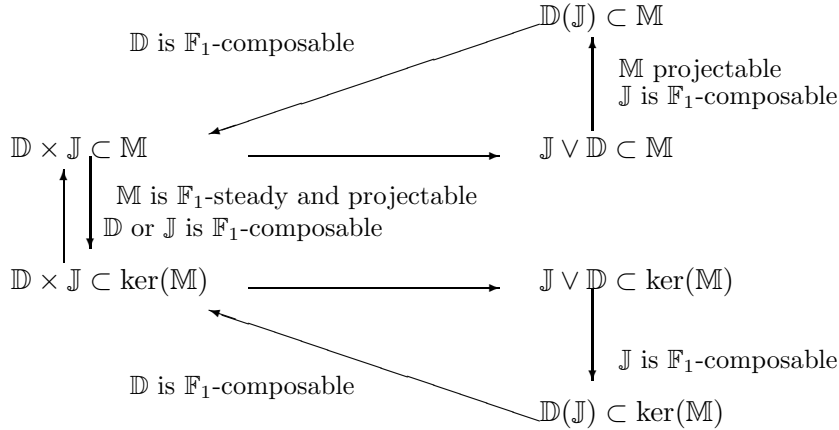
\begin{picture}(300,180)
\put(200,150){$\mathbb{D}(\mathbb{J}) \subset \mathbb{M}$}
\put(200,100){$\mathbb J \vee \mathbb D \subset \mathbb{M}$}
\put(0,100){$\mathbb D \times \mathbb J \subset \mathbb{M}$}
\put(0,50){$\mathbb D \times \mathbb J \subset \ker(\mathbb{M})$}
\put(200,50){$\mathbb J \vee \mathbb D \subset \ker(\mathbb{M})$}
\put(200,0){$\mathbb{D}(\mathbb{J}) \subset \ker(\mathbb{M})$}
\put(220,50){\vector(0,-1){35}}
\put(230,20){\small $\mathbb{J}$ is $\mathbb{F}_1$-composable}
\put(200,0){\vector(-3,1){125}}
\put(45,10){\small $\mathbb{D}$ is $\mathbb{F}_1$-composable}
\put(20,60){\vector(0,1){35}}
\put(30,100){\vector(0,-1){35}}
\put(90,50){\vector(1,0){95}}
\put(90,100){\vector(1,0){95}}
\put(220,110){\vector(0,1){35}}
\put(200,150){\vector(-3,-1){125}}
\put(35,70){\small $\mathbb{D}$ or $\mathbb{J}$ is $\mathbb{F}_1$-composable}
\put(40,82){\small $\mathbb{M}$ is  $\mathbb{F}_1$-steady and projectable}
\put(45,140){\small $\mathbb{D}$ is $\mathbb{F}_1$-composable}
\put(230,120){\small $\mathbb{J}$ is $\mathbb{F}_1$-composable}
\put(230,130){\small $\mathbb{M}$ projectable}
\end{picture}
\bigskip
\begin{cor}\label{cor:fulleq}
Let $\mathbb {J}$ and $\mathbb {D}$ be two $\mathbb
{F}_1$-composable classes of filters, and $\mathbb {M}$ be
an $\mathbb {F}_1$-steady projectable class of filters. The following
are equivalent:
\begin{enumerate}
\item 
$\mathbb {J}$ is $(\mathbb {D},\mathbb {M})$-composable;
\item 
$\mathbb {J}$ is $(\mathbb {D},\ker(\mathbb {M}))$-composable;
\item 
$\mathbb {J}$ is $(\mathbb {D},\ker(\mathbb {M}))$-steady;
\item 
$\mathbb {J}$ is $(\mathbb {D},\mathbb {M})$-steady;
\item 
$(\F\times\L)\in\ker(\mathbb {M})$ whenever $\L\in\mathbb {D}$ and
$\F\in\mathbb {J}$;
\item 
$(\F\times\L)\in\mathbb {M}$ whenever $\L\in\mathbb {D}$ and
$\F\in\mathbb {J}$.
\end{enumerate}
\end{cor}

In particular,  $\mathbb {F}_\omega$, $\mathbb
{F}_\omega^\bigtriangleup$, $\mathbb
{F}_\omega^{\bigtriangleup\bigtriangleup}$ are all $\mathbb
{F}_1$-composable, so that steady and composable are the same for
these classes. Therefore, Corollary \ref{cor:supstrFrechet}
combines with Corollary \ref{cor:fulleq} to the effect that

\begin{prop}\cite{quest}\cite{JM}
The following are equivalent:
\begin{enumerate}
\item $X$ is strongly Fr\'echet;
\item $X\times Y$ is strongly Fr\'echet for every productively Fr\'echet space $Y$;
\item $X\times Y$ is Fr\'echet for every first-countable space $Y$.
\end{enumerate}
\end{prop}
and
\begin{prop}\cite{JM}   \label{th:mainJM}
$X$ is productively Fr\'echet if and only if $X\times Y$ is Fr\'echet
(equivalently strongly Fr\'echet) for every strongly Fr\'echet space $Y$.
\end{prop}
This last result solves an old open question. See \cite{JM} for details.

Analogously, $\mathbb {F}_1^{\bigtriangleup}$ (the class of
Fr\'echet filters) is $\mathbb {F}_1$-composable.  From Theorem
\ref{th:main} applied with $\star=\bigtriangleup$ and
 $\mathbb {J}=\mathbb {F}_1$ we get
\begin{prop}\cite{myn.cont} \label{prop:fr}
A topological space is Fr\'echet if and only if its product with every
finitely generated space is Fr\'echet.
\end{prop}

The classes  $\mathbb {F}_1$, $\mathbb {F}_{\omega}$, and $\mathbb
F_\omega^\dagger$, and $\mathbb {F}_1^{\diamondsuit}$ are $\mathbb
{F}_1$-composable, so that steady and composable are the same for
these classes. Therefore, Corollary \ref{cor:suptight} and
Corollary \ref{cor:supfantight} combine with Corollary \ref{cor:fulleq}
to the effect that:

\begin{prop}\label{th:cow}
$X$ is countably tight if and only if $X\times Y$ is countably tight for every countably based $Y$.
\end{prop}
and
\begin{prop}\label{th:cow1}
$X$ is countably fan tight if and only if $X\times Y$ is countably fan tight for every countably based space $Y$.
\end{prop}

In contrast, we will see that the classes  $\mathbb
{F}_1^{\diamondsuit\diamondsuit}$, $\mathbb
{F}_1^{\bigtriangleup\diamondsuit}$, $\mathbb
{F}_1^{\dagger\diamondsuit}$, and $\mathbb
{F}_1^{\bigtriangleup\bigtriangleup}$  are not $\mathbb
{F}_1$-composable (\footnote{ We do not know if $\mathbb
{F}_{\omega}^{\bigtriangleup\dagger}$, $\mathbb
{F}_{\omega}^{\bigtriangleup\dagger}$, and $\mathbb
{F}_1^{\dagger\dagger}$ are $\mathbb {F}_1$-composable.}). Hence,
to investigate stability under product of the associated
properties, we need to introduce more machinery.

A $\vee$-compatible relation on filters (that is, defined on
$\mathbb {F}(X)$ for any set $X$) is $\times$-{\em compatible}
provided that
\begin{equation}\tag{$\bullet\bullet$}
(\F\times\G)\star A\then \G\star A\F.
\end{equation}
The proof of the following lemma is left to the reader.
\begin{lem}
$\bigtriangleup$, $\diamondsuit$ and $\dagger$ are all $\times$-compatible.
\end{lem}

\begin{lem}\label{lem:toto}
If $\mathbb {J}$ is an $\mathbb {F}_1$-composable class of filters
and $\star$ is a $\times$-compatible relation, then $\mathbb
{J}^{\star}$ is $\mathbb F_1$-steady and projectable.
\end{lem}
\begin{proof}
The fact that $\mathbb{J}^{\star}$ is $\mathbb F_1$-steady follows
from Proposition \ref{prop:stability}.
 Let $\F\in\mathbb
{J}^{\star}(X\times Y)$ and $\G\in\mathbb {J}(Y)$. Since
$X\times\G\in\mathbb {J}(X\times Y)$, $(X\times\G)\star\F$. Since
$\star$ is $\times$-compatible, $\F X\star\G$.  However, $\F
X=\pi_Y(\F)$.  Thus, $\pi_Y(\F)\in\mathbb {J}^{\star}$.
\end{proof}

\begin{thm} \label{th:kernel}
Let $\mathbb {J}$ be a steady $\mathbb {F}_1$-composable class of
filters containing $\mathbb {F}_1$,  and let $\star$, $\square$,
and $\bigtriangledown$ be $\times$-compatible relations.
\begin{enumerate}
\item The following are equivalent:
\begin{enumerate}\label{big1}
\item $\F\in\ker(\mathbb {J}^\star)$;
\item $(\F\times\G)\in\ker(\mathbb {J}^\star)$ for every $\G\in\ker(\mathbb {J}^{\star\star})$;
\item $(\F\times\G)\in\mathbb {J}^\star$ for every $\G\in\mathbb {J}$;
\item $(\F\times\G)\in\mathbb {F}_1^\star$ for every $\G\in\mathbb {J}$.
\end{enumerate}
\item If $\mathbb {J}^{\star}$ is $\mathbb {F}_1$-composable and contains $\mathbb {F}_1$, then the following are equivalent:
\begin{enumerate}\label{big2}
\item $\F\in\ker(\mathbb {J}^{\star\square})$;
\item $(\F\times\G)\in\ker((\ker(\mathbb {J}^{\star\star}))^\square)$ for every $\G\in\mathbb {J}^{\star}$;
\item $(\F\times\G)\in\mathbb {J}^\square$ for every $\G\in\mathbb {J}^{\star}$;
\item $(\F\times\G)\in\mathbb {F}_1^\square$ for every $\G\in\mathbb {J}^{\star}$.
\end{enumerate}
\item \label{big3}
Let $\mathbb {J}^{\star}$ be an $\mathbb {F}_1$-composable class of filters.  If $\F\in \ker(\mathbb {J}^{\bigtriangledown\square})$ then $(\F\times\G)\in\ker(\mathbb {J}^{\star\square})$ for every $\G\in\ker(\mathbb {J}^{\star\bigtriangledown})$.
\end{enumerate}
\end{thm}
\begin{proof}
We prove (\ref{big1}).

(1a) $\then$ (1b).  By Theorem
\ref{th:main} (1), $\F\vee\G\in \mathbb J^\star$ for every $\G\in \mathbb J^{\star\star}$, hence in particular for every $\G\in\ker(\mathbb J^{\star\star})$. Since $\mathbb J$ is $\mathbb F_1$-composable, Lemma \ref{lem:toto} applies to the effect that $\mathbb J^\star$ is $\mathbb F_1$-steady and projectable. Moreover, $\ker (\mathbb J^\star)$ and $\ker(\mathbb J^{\star\star})$ are $\mathbb F_1$-composable by definition. In view of Corollary \ref{cor:fulleq}, (1b) follows.

(1b) $\then $ (1c) because $\mathbb J\inc\ker(\mathbb J^{\star\star})$ and $\ker(\mathbb J^{\star})\inc \mathbb J^\star$. Similarly,
(1c) $\then$ (1d) because $\mathbb J^\star\inc \mathbb F_1^\star$.

(1d) $\then$ (1a). Let $A\inc X\times Y$. We want to show that $A\F\in \mathbb J^\star$. If $\G\in \mathbb J(Y)$ meshes with $A\F$ then $A\mesh(\F\times\G)$. But $\F\times\G\in \mathbb F_1^\star$, so that $(\F\times\G)\star A$. By $(\bullet\bullet)$, we have $\G\star A\F$, so that $A\F\in \mathbb J^\star$.

 (2a) $\then$ (2b) is proved from Theorem \ref{th:main} (2), combined with Lemma \ref{lem:toto} and Corollary \ref{cor:fulleq} as (1a) $\then$ (1b). (2b) $\then$ (2c) $\then$ (2d) follows from inclusions among the classes of filters considered and (2d) $\then$ (2a) follows from $(\bullet\bullet)$.

Finally (\ref{big3}) follows in a similar way from Corollary \ref{cor:comp1}.
\end{proof}
In the very same way that we deduced Corollary \ref{cor:sup} from Theorem \ref{th:main}, we deduce the following from Theorem \ref{th:kernel}.
\begin{cor} \label{cor:product}
Let $\mathbb {J}$ be a steady $\mathbb {F}_1$-composable class containing
$\mathbb {F}_1$ and let $\star$ and $\square$ be two
$\times$-compatible relations.
\begin{enumerate}
\item The following are equivalent:
\begin{enumerate}
\item
$X$ is $\ker(\mathbb {J}^\star)$-based;
\item $X\times Y$ is $\ker(\mathbb {J}^\star)$-based for every $\ker(\mathbb
{J}^{\star\star})$-based $Y$;
\item $X\times Y$ is $\mathbb {J}^\star$-based for every $\mathbb
{J}$-based $Y$;
\item $X\times Y$ is $\mathbb {F}_1^\star$-based for every $\mathbb
{J}$-based $Y$.
\end{enumerate}
\item If $\mathbb {J}^{\star}$ is $\mathbb {F}_1$-composable and contains $\mathbb {F}_1$, then the following are equivalent:
\begin{enumerate}
\item $X$ is $\ker(\mathbb {J}^{\star\square})$-based;
\item $X\times Y$ is $\ker((\ker(\mathbb {J}^{\star\star}))^\square)$-based for every $\mathbb {J}^{\star}$-based $Y$;
\item $X\times Y$ is $\mathbb {J}^\square$-based for every $\mathbb {J}^{\star}$-based $Y$;
\item $X\times Y$ is $\mathbb {F}_1^\square$-based for every $\mathbb {J}^{\star}$-based $Y$.
\end{enumerate}
\item If $\mathbb {J}^{\star}$ is $\mathbb {F}_1$-composable and contains $\mathbb {F}_1$, and if $X$ is $\ker(\mathbb {J}^{\bigtriangledown\square})$-based, then
 $X\times Y$ is $\ker(\mathbb {J}^{\star\square})$-based for every $\ker(\mathbb {J}^{\star\bigtriangledown})$-based $Y$.
\end{enumerate}
\end{cor}


\section{Characterizations of kernels}

Now the missing point to apply the previous results to theorems of stability of local
topological properties like tightness or fan-tightness under product is to characterize
kernels of the corresponding classes of filters.

\begin{thm} \label{th:kerprincipal}
$$\ker(\mathbb {F}_1^{\bigtriangleup\bigtriangleup})=\mathbb {F}_1.$$
\end{thm}
\begin{proof}
$\mathbb {F}_1\subseteq \ker(\mathbb
{F}_1^{\bigtriangleup\bigtriangleup})$ is clear. Conversely, if
$X$ is a set and $\F\notin \mathbb {F}_1(X)$, then for every
$F\in\F$, there exists $H_F\in\F$ and
 $x_F\in F\setminus H_F$. Let $Y$ be an infinite set. Then $\F\times Y$ is not an almost principal filter because the sets $\{x_F\}\times Y$ are infinite. By Theorem \ref{th:almostprinc},
$\F\times Y\notin \mathbb {F}_1^{\bigtriangleup\bigtriangleup}$, and $\F\notin \ker(\mathbb F_1^{\bigtriangleup\bigtriangleup})$.
\end{proof}

As a consequence, we obtain the dual statement to Proposition \ref{prop:fr}.
\begin{cor}\cite{myn.cont}
A topological space is finitely generated if and only if its product with every
Fr\'echet space is Fr\'echet.
\end{cor}
\subsection{$\kappa$-tightness, productive $\kappa$-tightness, tight points, and absolute countable tightness}
For our discussion of kernels for tightness we consider the
following relation. Let $\kappa$ be an infinite cardinal and
define the relation $\diamondsuit_{\kappa}$ on $\mathbb {F}(X)$ by
$\F{\diamondsuit_{\kappa}}\H$ if
$$
\H\mesh\F \then \exists A \in \mathbb {F}_1, |A|\leq\kappa: A\mesh \F\vee\H.
$$
This is a $\vee$-compatible and $\times$-compatible relation. 
$\mathbb {F}_1^{\diamondsuit_{\kappa}}$ is the class of {\em
$\kappa$-tight filters}. Recall that a topological space $X$ has
tightness $\kappa$ \cite{Ar.spectrum} if for every $x\in X$ and
$A\subseteq X$ such that $x\in \cl A$,
 there exists a subset
$B\subseteq A$ such that $|B|\leq\kappa$ and  $x\in \cl B$. It is easy to see that a space $X$ has $\kappa$-tightness if and only if all its neighborhood filters are $\kappa$-tight.

Now we give a general result for kernels of
$\diamondsuit_{\kappa}$-polars for classes of filters satisfying certain conditions. 

Let $I$ and $J$ be sets.  A function $\gamma:I\times J\to 2^X$ is called a {\em presentation on} $X$.  
We define 
$$\gamma_*=\{\bigcup_{\alpha\in I} \gamma(\alpha,\beta): 
\beta\in J\}^\uparrow$$
If $\G$ is a filter on $X$ and 
$\G=\gamma_*$, then we say that $\gamma$ is a {\em presentation of} $\G$.  Of course, if $\gamma_*$ is a filter, then $\gamma$ is a presentation of $\gamma_*$. 
A presentation $\gamma:I\times J\to 2^X$ is called $\mathbb J$-{\em expandable} if the filter 
$\gamma^*$ defined on $X\times I$ by
 $$\gamma^*=\{ \{\bigcup_{\alpha\in I} 
\gamma(\alpha,\beta)\times\{\alpha\}\}:\beta\in J\}^\uparrow $$
  is a $\mathbb J$-filter on $X\times I$. Notice that $\gamma_*=\pi_X(\gamma^*)$. Therefore, if $\mathbb J$ is projectable and $\gamma:I\times J\to 2^X$ is a $\mathbb J$-expandable presentation on $X$, then $\gamma_*\in \mathbb J(X)$.  

\begin{thm}\label{thm:generaltightpolar}
Let $\mathbb J$ be an $\mathbb F_1$-composable class of filters included in $\mathbb F_1^{\diamondsuit_\kappa}$. The following are equivalent:
\begin{enumerate}
\item \label{it:ker}
$\F\in\ker(\mathbb J^{\diamondsuit_\kappa});$
\item \label{it45}
$\F\in\mathbb {J}^{\diamondsuit_{\kappa}}$ and
$A\vee\F\in\ker(\mathbb {J}^{\diamondsuit_{\kappa}})$ for all $A$
such that $A\mesh\F$ and $|A|\leq\kappa$;
\item \label{it:presentation1}
For any $\mathbb J$-expandable presentation  $\gamma:I\times J\to 2^X$ such that $\gamma_*\mesh\F$, there exists a subset $C$ of $I$ of cardinality at most $\kappa$ such that 
$\{\bigcup_{\alpha\in C} \gamma(\alpha,\beta):\beta\in 
J\}\mesh\F$;
\item \label{it:presentation}
For any $\mathbb J$-expandable presentation  $\gamma:\F\times \J\to 2^X$  such that $\gamma_*\mesh\F$, where $\J$ is a $\mathbb J$-filter and each set $\gamma(F,J)$ is a subset of $F$ of cardinality at most $\kappa$, there exists a subfamily $\K$ of $\F$ of cardinality at most $\kappa$ such that 
$\{\bigcup_{F\in\K} \gamma(F,J):J\in \J\}\mesh\F$.
\end{enumerate}
\end{thm}

\begin{proof}
(\ref{it:ker}) $\then$ (\ref{it45}) is obvious.

(\ref{it45}) $\then$ (\ref{it:presentation1}). 
Let $\gamma:I\times J\to 2^X$ be a $\mathbb J$-expandable presentation such that $\gamma_*\mesh\F$. As $\F\in\mathbb J^{\diamondsuit_\kappa}$, there exists $A$ of cardinality at most $\kappa$ such that $A\mesh\F\vee\gamma_*$ or equivalently, $\F\vee A\mesh\gamma_*$. By $\mathbb J$-expandability of $\gamma$, the filter $\gamma^*$ is a $\mathbb J$-filter on $X\times I$. 
Moreover $\pi_X^-(\F\vee A)\in \mathbb J^{\diamondsuit_\kappa}(X\times I)$ because $\F\vee A\in\ker(\mathbb {J}^{\diamondsuit_{\kappa}})$; and $\pi_X^-(\F\vee A)\mesh\gamma^*$. Thus, there exists a set $D$ of cardinality at most $\kappa$ such that $D\mesh\pi_X^-(\F\vee A)\vee\gamma^*$. 
Then $\{\bigcup_{\alpha\in \pi_I(D)} \gamma(\alpha,\beta):\beta\in J\}\mesh\F$.

(\ref{it:presentation1}) $\then$ (\ref{it:presentation}) is obvious.

(\ref{it:presentation}) $\then$ (\ref{it:ker}).
Let $A\inc X\times Y$ and let $\G$ be a $\mathbb J$-filter on $Y$ such that $\G\mesh A\F$. For every $F\in\F$, $\G\mesh AF$ and $\G\in \mathbb F_1^{\diamondsuit_\kappa}$ so that there exists a subset $B_F$ of $AF$ of cardinality at most $\kappa$ and meshing with $\G$. For each $F\in \F$, there exists a function $f_F:B_F\to F$ with graph included in $A$. 
Consider the function $\gamma:\F\times\G\to 2^X$ defined by $\gamma(F,G)=f_F(G\cap B_F)$.  It is a $\mathbb J$-expandable presentation, because $\mathbb J$ is $\mathbb F_1$-composable and
$\{\bigcup_{F\in\F} \gamma(F,G)\times\{F\}: G\in\G\}^\uparrow$ is 
the image of the $\mathbb J$-filter $\G$ under the multivalued map $R:Y\rightrightarrows X\times \F$ with graph $\{(y,x,F):f_F(x)=y\}$.  By (\ref{it:presentation}), there is a subfamily $\K$ of $\F$ of cardinality at most $\kappa$ such that 
$\{\bigcup_{F\in \K} \gamma(F,G):G\in \G\}\mesh\F$. The set 
$C=\bigcup_{F\in\K} B_F$ is of cardinality at most $\kappa$ 
and $C\mesh\G\vee A\F$.

\end{proof}

Let $\bigwedge(\mathbb{D})$ be the class of filters than can be 
represented as the infimum of a family of $\mathbb{D}$-filters.
\begin{cor}\label{cor:tightgeneral}
If $\mathbb{J}$ is an $\mathbb{F}_1$-composable class included
in the class $\mathbb{F}_1^{\diamondsuit_{\kappa}}$ such that 
$\mathbb{J} = \bigwedge(\mathbb{D})$ where $\mathbb{D}$ contains 
$\mathbb{F}_1$, 
then the following are equivalent:
\begin{enumerate}
\item \label{itcor:ker}
$\F\in\ker(\mathbb{J}^{\diamondsuit_\kappa});$
\item \label{it45b}
$\F\in\mathbb{J}^{\diamondsuit_{\kappa}}$ and
$A\vee\F\in\ker(\mathbb{J}^{\diamondsuit_{\kappa}})$ for all $A$
such that $A\mesh\F$ and $|A|\leq\kappa$;
\item \label{it41b}
for every family $(\G_\alpha)_{\alpha\in I}$ of 
$\mathbb{D}$-filters such that $\bigwedge_{\alpha\in I} 
\G_{\alpha}\mesh\F$ there
exists  $J\subseteq I$  such that $|J|\leq\kappa$ and
$\bigwedge_{\alpha\in J} \G_{\alpha}\mesh\F$;
\item \label{it42b}
for every family $(\G_\alpha)_{\alpha \in I}$ of $\mathbb{D}$-filters: if 
$\forall F\in \F$, $\exists \alpha \in I$ and 
$C_F\subseteq F$: $C_F\in \G_\alpha$ and $|C_F|\leq \kappa$, then 
there exists $J\subseteq I$ such that $|J|\leq\kappa$ and
$\bigwedge_{\alpha\in J} \G_{\alpha}\mesh\F$;
\item \label{it43b}
for every family $(\G_\alpha)_{\alpha \in I}$ of $\mathbb{D}$-filters, 
each of which is either free or the principal filter of a set of 
cardinality at most $\kappa$: if $\forall F\in \F$,
$\exists \alpha \in I$  and $C_F\subseteq F$: 
$C_F\in \G_\alpha$ and $|C_F|\leq \kappa,$ then there exists 
$J\subseteq I$ such that $|J|\leq\kappa$ and 
$\bigwedge_{\alpha\in J} \G_{\alpha}\mesh\F$.
\end{enumerate}
\end{cor}
\begin{proof}
(\ref{itcor:ker}) $\iff$ (\ref{it45b}) follows from Theorem \ref{thm:generaltightpolar} and (\ref{it41b}) $\then$ (\ref{it42b}) $\then$ (\ref{it43b}) by definition.

(\ref{itcor:ker}) $\then$ (\ref{it41b}).
Let $(\G_\alpha)_{\alpha\in I}$ be a family of $\mathbb
{D}$-filters such that $\bigwedge_{\alpha\in I} \G_{\alpha}\mesh\F$. 
Consider the presentation 
$\gamma:I\times \prod_{\alpha\in I} \G_\alpha\to 2^X$ defined by 
$\gamma(\alpha,(G_\beta)_{\beta\in I})=G_\alpha$. It is a 
$\mathbb{J}$-expandable presentation of 
$\bigwedge_{\alpha\in I} \G_\alpha$.  
Indeed, $\mathbb J$ is $\mathbb F_1$-composable so that $\G_{\alpha}\times\{\alpha\}\in\mathbb J$ for all $\alpha\in I$, and $\mathbb J$ stable by infimum so that 
$\bigwedge_{\alpha\in I} (\G_\alpha \times\{\alpha\})$ is a 
$\mathbb J$-filter on $X\times I$. In view of Theorem \ref{thm:generaltightpolar}, there exists $C\inc I$ such that $|C|\leq\kappa$ and 
$\bigwedge_{\alpha\in C} \G_\alpha\mesh\F$.

(\ref{it43b}) $\then$ (\ref{itcor:ker}).
In view of Theorem \ref{thm:generaltightpolar}, we only need to show that for any $\mathbb J$-filter $\G$ such that $\G\mesh\F$ and every $\mathbb J$-expandable presentation  $\gamma:\F\times \J\to 2^X$ of $\G$ such that $\J$ is a $\mathbb J$-filter and each $\gamma(F,J)$ is a subset of $F$ of cardinality at most $\kappa$, there exists a subfamily $\K$ of $\F$ of cardinality at most $\kappa$ such that 
$\{\bigcup_{F\in \K} \gamma(F,J):J\in \J\}\mesh\F$.
For such a presentation of $\G$, define for each $F\in\F$ the filter $\J_F=\gamma(F,\J)$. By $\mathbb F_1$-composability of $\mathbb J$, it is a $\mathbb J$-filter. As $\mathbb J=\bigwedge(\mathbb D)$, there exists a $\mathbb D$-filter $\L_F\geq\J_F$ that is either free
or principal. Indeed, if $\J_F^\bullet$ is non degenerate, then it is a $\mathbb D$-filter finer than $\J_F$. Otherwise, $\J_F$ is free and therefore admits finer free $\mathbb D$-filters. 
$\bigwedge_{F\in\F} {\L_F}\mesh\F$ and each $\L_F$ 
contains a subsets a cardinality at most $\kappa$ of $F$. Therefore, there 
exists a subfamily $\K$ of $\F$ of cardinality at most $\kappa$ such that 
$\bigwedge_{F\in\K} {\L_F}\mesh\F$. Moreover, 
$\{\bigcup_{F\in\K} \gamma(F,J):J\in\J\}^\uparrow
\leq \bigwedge_{F\in\K} {\J_F}
\leq \bigwedge_{F\in\K} {\L_F}$. Therefore, 
$\{\bigcup_{F\in\K} \gamma(F,J):J\in\J\}^\uparrow\mesh\F$.
\end{proof}

\subsubsection{Productive $\kappa$-tightness}

\begin{cor} \label{cor:tightkernel}
The following are equivalent
\begin{enumerate}
\item\label{it1}
$\F\in \ker(\mathbb{F}_1^{\diamondsuit_{\kappa}\diamondsuit_{\kappa}})$;
\item\label{it0}
$\F\in\mathbb {F}^{\diamondsuit_{\kappa}\diamondsuit_{\kappa}}$ and 
$A\vee\F \in \ker(\mathbb{F}_1^{\diamondsuit_{\kappa}\diamondsuit_{\kappa}})$ 
for all $A$ such that $A\mesh\F$ and $|A|\leq\kappa$;
\item\label{it2}
for every collection 
$\{\G_\alpha\colon\alpha \in I\}
\subseteq \mathbb{F}_1^{\diamondsuit_{\kappa}}$: 
if $\bigwedge_{\alpha \in I} \G_\alpha\mesh\F$, then 
$\exists J\subseteq I$ such that $|J|\leq\kappa$ and 
$\bigwedge_{\alpha\in\J} \G_{\alpha}\mesh\F$;
\item\label{it3}
for every collection 
$\{\G_\alpha\colon\alpha \in I\}
\subseteq \mathbb{F}_1^{\diamondsuit_{\kappa}}$: 
if $\forall F\in \F$, $\exists\alpha \in I$ and $C_F\subseteq F$ such that 
$C_F\in \G_\alpha$ and $|C_F|\leq \kappa$, then $\exists J\subseteq I$ 
such that $ |J|\leq\kappa$ and 
$\bigwedge_{\alpha\in J} \G_{\alpha}\mesh\F$;
\item\label{it4}
for every collection $\{\G_\alpha\colon\alpha \in I\}
\subseteq \mathbb{F}_1^{\diamondsuit_{\kappa}}$ 
each of which is either free or the principal filter of a set of 
cardinality at most $\kappa$: 
if $\forall F\in \F$, $\exists \alpha \in I$ and 
$C_F\subseteq F$  such that $C_F\in \G_\alpha$
$|C_F|\leq \kappa$, then there exists
$J\subseteq I$ such that $|J|\leq\kappa$ and
$\bigwedge_{\alpha\in J} \G_{\alpha}\mesh\F$.
\end{enumerate}
\end{cor}

\begin{proof}
Since $\mathbb
{F}_1^{\diamondsuit_{\kappa}}=\bigwedge\left(\mathbb
{F}_1^{\diamondsuit_{\kappa}}\right)$, the corollary is just a
restatement of Corollary~\ref{cor:tightgeneral} with $\mathbb
{D}=\mathbb {J}=\mathbb {F}_1^{\diamondsuit_{\kappa}}$.
\end{proof}

A point $x$ of a topological space $X$ is {\em a  productively} $\kappa$-{\em tight point} if $\N_X(x)\in \ker(\mathbb F_1^{\diamondsuit_\kappa\diamondsuit_\kappa})$. A topological space is {\em   productively} $\kappa$-{\em tight} if all its points are productively $\kappa$-tight.

The interest of such an explicit description of a kernel lies in its combination with 
the corresponding instance of Corollary \ref{cor:product}. In the present case,
\begin{cor} \label{cor:tight}
$X\times Y$ is $\kappa$-tight for every $\kappa$-tight space $Y$ if and only if
 $X$  is productively $\kappa$-tight.
\end{cor}

Notice that the fourth condition in Corollary \ref{cor:tightkernel} corresponds to the non existence of
{\em singular families}  in the sense of Arhangel'skii \cite{Ar.spectrum}.
Hence, the combination of Corollary~\ref{cor:tightkernel}
 and Corollary \ref{cor:tight} extends the equivalence
between b) and c) in \cite[Theorem 3.6]{Ar.spectrum}, and provides a shorter proof.  Of all the conditions of Corollary~\ref{cor:tightkernel} the third is probably the most straightforward and easy to use, but was apparently not known to be equivalent to productive $\kappa$-tightness.

\subsubsection{Tight points}

We now characterize $\ker(\mathbb
{F}_1^{\bigtriangleup\diamondsuit})$. For a cardinal $\kappa$ we
let $S_{\kappa}$ denote the space obtained as a quotient of
$\kappa$-many mutually disjoint convergent sequences
 $(z^{\xi}_n)_{n\in\omega}$ by identifying their limit points
 to a single point $w$. $S_\kappa$ is endowed with the quotient topology.
  We denote the neighborhood of $w$ by $\S_{\kappa}$.
\begin{thm}\label{thm:kerbigtrwr}
The following are equivalent for a filter on $X$:
\begin{enumerate}
\item\label{it31}
$\F\in\ker(\mathbb {F}_1^{\bigtriangleup\diamondsuit})$;
\item\label{it34}
$\F\in\mathbb {F}_1^{\bigtriangleup\diamondsuit}$ and $\F\vee
A\in\ker(\mathbb {F}_1^{\bigtriangleup\diamondsuit})$ for any
countable set $A$ meshing with $\F$;
\item\label{it30}
If $\bigwedge_{\alpha\in I}(x_n^{\alpha})_{n\in\omega}\mesh\F$,
 then there exists a countable $J\subseteq I$ such that
$\bigwedge_{\alpha\in J}(x_n^{\alpha})_{n\in\omega}\mesh\F$;
\item\label{it35}
If $\bigwedge_{\alpha\in I}(x_n^{\alpha})_{n\in\omega}\mesh\F$,
 where each sequence $(x_n^{\alpha})_{n\in\omega}$ is either free or principal, 
 then there exists a countable $J\subseteq I$ such that
$\bigwedge_{\alpha\in J}(x_n^{\alpha})_{n\in\omega}\mesh\F$;

\item\label{it32}
If $\{\G_{\alpha}\colon\alpha\in I\}\subseteq\mathbb
{F}_1^{\bigtriangleup}$ and  $\bigwedge_{\alpha\in
I}\G_{\alpha}\mesh\F$, then there exists a countable $J\subseteq I$
such that $\bigwedge_{\alpha\in J}\G_{\alpha}\mesh\F$;
\item\label{it33}
$\F\times\S_{|X^{\omega}|}\in\mathbb {F}_1^{\diamondsuit}$.
\end{enumerate}
\end{thm}
\begin{proof}
(\ref{it31}) $\iff$(\ref{it30}) $\iff$ (\ref{it35})   $\iff$ (\ref{it34}) follows from
Corollary~\ref{cor:tightgeneral} with $\mathbb {J}=\mathbb
{F}_1^\bigtriangleup$ and $\mathbb {D}$ equal to the class of
sequences and the observation that $\mathbb
{F}_1^\bigtriangleup=\bigwedge(\mathbb {D})$ (e.g., \cite{dolseq}).

(\ref{it31}) $\iff$ (\ref{it32}) follows from
Corollary~\ref{cor:tightgeneral} with $\mathbb {J}=\mathbb
{F}_1^\bigtriangleup$ and $\mathbb {D} =\mathbb
{F}_1^\bigtriangleup$ and the observation that $\mathbb
{F}_1^\bigtriangleup=\bigwedge(\mathbb {F}_1^\bigtriangleup)$.

(\ref{it31}) $\then$ (\ref{it33}). By Theorem~\ref{th:eq},
$\F\times\G\in\mathbb {F}_1^{\diamondsuit}(X\times Y)$
 for any Fr\'echet filter $\G$ on any set $Y$.  In particular,
 $\F\times\S_{|X^{\omega}|}\in\mathbb {F}_1^{\diamondsuit}(X\times S_{|X^{\omega}|})$.

(\ref{it33}) $\then$ (\ref{it30}).
Let  $\{(x^{\alpha}_n)_{n\in\omega}\colon \alpha\in I\}$ be a collection of
 sequences such that $(\bigwedge_{\alpha\in I}(x_n^{\alpha})_{n\in\omega})\mesh\F$.
 Let $\kappa=|X^{\omega}|$ and $K=\{(x^{\xi}_{n})_{n\in\omega}\colon\xi\in\kappa\}$
 be an enumeration of the set $\{(x^{\alpha}_n)_{n\in\omega}\colon \alpha\in I\}$
possibly with repetitions.  Let $R\subseteq X\times S_{\kappa}$ be defined by
 $R=\bigcup_{\xi\in\kappa}\{(x^{\xi}_n,z^{\xi}_n)\colon n\in\omega\}$.
 Notice that $R\mesh(\F\times\S_{\kappa})$.  By (\ref{it33}),
 there is a countable set $T\mesh(R\vee(\F\times\S_{\kappa}))$.
 Let $Q=\{\xi\in\kappa\colon(\exists n\in\omega)(x^{\xi}_n,z^{\xi}_n)\in T\}$.
Let $F\in\F$.
Since $T\mesh(R\vee(F\times\S_{\kappa}))$, there is  $\xi\in Q$ such that
 $\pi_X^{-1}(F)\mesh(x^{\xi}_n,z^{\xi}_n)_{n\in\omega}$.  So, for some $\xi\in Q$
we have $F\mesh(x^{\xi}_n)_{n\in\omega}$.  Thus,
$\F\mesh\bigwedge_{\xi\in Q}(x^{\xi}_n)_{n\in\omega}$ and $Q$ is countable.
\end{proof}

In \cite{bella.vanmill} the notion of a tight point is introduced.  A
 collection of sets $\E$ is said to {\em cluster at a point} $p$ provided that for
 every neighborhood $U$ of $x$ there is an $E\in\E$ such that $U\cap E$ is infinite.
  A point $x$ in a space $X$ is said to be {\em tight} provided that for any
collection of sets $\E$ that clusters at $p$ one can find a countable
subcollection
$\E^*\subseteq\E$ such that $\E^*$ clusters at $p$.

\begin{thm}\label{thm:tight=ker}
A point $p$ of a space $X$ is tight if and only if $\N_X(p)\in\ker(\mathbb {F}_1^{\bigtriangleup\diamondsuit})$.
\end{thm}

\begin{proof}
Suppose that the neighborhood filter $\N(p)$ is in $\ker(\mathbb
{F}_1^{\bigtriangleup\diamondsuit})$. Let $\E$ be a collection of
sets which clusters at $p$. For each $F\in\N(p)$ we can find an
$E^F\in\E$ such that $E^F\cap F$ is infinite.  Let
$(x^F_n)_{n\in\omega}$ be a
 free sequence on $E^F\cap F$.  Clearly,
$\N(p)\mesh\bigwedge_{F\in\N(x)}(x^{F}_n)_{n\in\omega}$.
There exists $\{F_k\colon k\in\omega\}\subseteq \F$ such that
$\N(p)\mesh\bigwedge_{k\in\omega}(x^{F_k}_n)_{n\in\omega}$.
Let $\E^*=\{E^{F_k}\colon k\in\omega\}$.
 It is easily verified that $\E^*$ clusters at $p$.

Suppose that $p$ is a tight point of $X$.
Let $\N(p)\mesh\bigwedge_{\alpha\in I}(x^{\alpha}_n)_{n\in\omega}$, where each sequence $(x^{\alpha}_n)_{n\in\omega}$ is either free
of principal. In other words, we can write  
$\N(p)\mesh A \bigwedge_{\alpha\in J} (x^{\alpha}_n)_{n\in\omega}$, 
where  $A=\bigwedge\{(x^{\alpha}_n)_{n\in\omega}\colon(x^{\alpha}_n)_{n\in\omega} \in\mathbb F_1\}$ and every sequence $(x^{\alpha}_n)_{n\in\omega}$ with $\alpha\in J$ is free. 
If $A\mesh\N(p)$, either  $A\mesh\N(p)^\bullet$ or $A\mesh\N(p)^\circ$. In the first case, there is a fixed sequence of the original collection whose defining point is in $\N(p)^\bullet$ and is therefore finer than $\N(p)$. In the second case, $A\cap N$ is infinite for every $N\in\N(p)$ and there exists a countably infinite set  $A_N\inc A\cap N$. The collection $\{A_N:N\in\N(p)\}$ clusters at $p$. Since $p$ is a tight point, there exists $\{A_{N_k}:k\in\omega\}$ which clusters at $p$. 
Then $\bigcup_{k\in\omega} A_{N_k}$ is a countable subset 
of $A$ 
meshing with $\N(p)$, and therefore defines a countable collection of (fixed) sequences from the original collection 
whose infimum is meshing with $\N(p)$.

If $A(\neg\mesh)\N(p)$ then $\N(p)\mesh\bigwedge_{\alpha\in J}(x^{\alpha}_n)_{n\in\omega}$. Each $N\in\N(p)$ is meshing with one of the sequences 
$(x^{\alpha}_n)_{n\in\omega}$ for $\alpha\in J$, so that the family $\{\{x^{\alpha}_n: n\in\omega\}:\alpha\in J\}$ clusters at $p$. Since $p$ is tight, there is $(\alpha_k)_{k\in\omega}$ in $J$ such that $\{\{x^{\alpha_k}_n:n\in\omega\}:k\in\omega\}$ clusters at $p$. Then 
$\N(p)\mesh\bigwedge_{k\in\omega}(x^{\alpha_k}_n)_{n\in\omega}$.
\end{proof}

From Theorem~\ref{thm:kerbigtrwr} and Theorem~\ref{thm:tight=ker}, we
 have the following improvement of  \cite[Theorem 3.1 ]{bella.vanmill}
which states the equivalence between (1) and (2), but only for a {\em countable} space $X$.

\begin{cor}\label{cor:card}
The following are equivalent:
\begin{enumerate}
\item every point of $X$ is tight;
\item $X\times S_{|X^{\omega}|}$ is countably tight;
\item $X\times Y$ is countably tight for every Fr\'echet space $Y$.
\end{enumerate}
\end{cor}

By \cite[Example 1.3]{Ar.spectrum}, $S_{\omega}\times
S_{\omega^{\omega}}$ is not countably tight.  Since $S_{\omega}$
is a countable space, $\S_{\omega}\in\mathbb
{F}_1^{\diamondsuit\diamondsuit}$.  As
$\S_{\omega^{\omega}}\in\mathbb {F}_1^{\bigtriangleup}$, we have
$\mathbb {F}_1^{\diamondsuit\diamondsuit}\setminus\ker(\mathbb
{F}_1^{\bigtriangleup\diamondsuit})\neq\emptyset$. Hence, $\mathbb
{F}_1^{\diamondsuit\diamondsuit}\neq\ker(\mathbb
{F}_1^{\diamondsuit\diamondsuit})$ and $\mathbb
{F}_1^{\bigtriangleup\diamondsuit}\neq\ker(\mathbb
{F}_1^{\bigtriangleup\diamondsuit})$.

In view of Theorem \ref{thm:tight=ker} and considering that $\mathbb F_1^{\bigtriangleup} \inc \mathbb F_1^{\diamondsuit}$, hence
$\ker(\mathbb {F}_1^{\diamondsuit\diamondsuit})\inc \ker(\mathbb {F}_1^{\bigtriangleup\diamondsuit})$, we obtain:

\begin{cor}\label{cor:prodtighttotight}
Every productively countably tight point is a tight point.
\end{cor}
\subsubsection{Absolute tightness}
A notion related to tightness, productive tightness and tight points and introduced in \cite{Ar.spectrum} is that of absolute tightness. A point of a topological space $X$ is {\em absolutely tight} if it is a point of countable tightness in a compactification of $X$. Absolute tightness can be characterized in a way that is similar to the characterization of productive tightness in Corollary \ref{cor:tightkernel}.

\begin{thm}\label{thm:abstight}
Let $X$ be a completely regular topological space and let $x\in X$. The following are equivalent:
\begin{enumerate}
\item $x$ is an absolutely tight point;
\item for any family $\{\F_\alpha:\alpha\in I\}$ of filters such that 
$\bigwedge_{\alpha\in I} \F_\alpha\mesh\N(x)$ there 
exists 
$(\alpha_i)_{i\in\omega}$ in $I$ such that 
$\bigwedge_{i\in\omega} \F_{\alpha_i}\mesh\N(x)$;
\item for any family $\{\U_\alpha:\alpha\in I\}$ of ultrafilters such that 
$\bigwedge_{\alpha\in I} \U_\alpha\mesh\N(x)$ there exists 
$(\alpha_i)_{i\in\omega}$ in $I$ such that 
$\bigwedge_{i\in\omega} \U_{\alpha_i}\mesh\N(x)$.
\end{enumerate}
\end{thm}
\begin{proof}
(1) $\then$ (2).

Let $\{\F_\alpha:\alpha\in I\}$ be a family of filters such that 
$\bigwedge_{\alpha\in I} \F_\alpha\mesh\N(x)$. Since each 
$\F_\alpha = \bigwedge_{\U\in\beta(\F_\alpha)} \U$, where 
$\beta(\F)$ denotes the set of ultrafilters finer than $\F$, we have 
$\bigwedge_{\alpha\in I,\, \U\in\beta(\F_\alpha)} \U\mesh\N(x)$. 
Let $bX$ denote a compactification of $X$ and let
$A=\bigcup_{\alpha\in I,\, \U\in\beta(\F_\alpha)} \lm_{bX}\U$. 
Let  
$W\in\N_{bX}(x)$. By regularity of $bX$, there exists $V\in\N_{bX}(x)$ such that $\cl_{bX}V\inc W$. Since $V\cap X\in\N_X(x)$, there exists $\alpha\in I$, $\U\in\beta(\F)$ and $U\in\U$ such that $U\inc V\cap X$. Then $\lm_{bX}\U\inc \cl_{bX}V\inc W$. Therefore, $x\in \cl_{bX} A$. Since $x$ is an absolutely tight point, there exist points $a_n$ in $A$ such that $x\in\cl_{bX}(\{a_n:n\in\omega\})$. For each $n$, pick $\U_n$ in 
$\bigcup_{\alpha\in I} \beta(\F_\alpha)$ such that 
$a_n\in\lm_{bX}\U_n$, and pick $\alpha_n$ such that $\U_n\in\beta(\F_{\alpha_n})$. We claim that 
$\N(x)\mesh \bigwedge_{n\in\omega} \F_{\alpha_n}$. 
Indeed, for each open $B\in\N_X(x)$, there is an open $B_1\in\N_{bX}(x)$ such that $ B_1\cap X=B$. There is an $n\in\omega$ such that $a_n\in B_1$. Therefore, there exists $U\in \U_n$ such that $U\inc B_1\cap X=B$. Hence 
$\N(x) \mesh \bigwedge_{n\in\omega} \U_{_n}$ so that 
$\N(x) \mesh \bigwedge_{n\in\omega} \F_{\alpha_n}$. 

(2) $\then$ (3) is obvious.

(3) $\then$ (1).
Let $bX$ denote a compactification of $X$ and let $A\inc bX$ be such that $x\in \cl_{bX}A$. Let $W\in\N_X(x)$ be open. There is a $bX$-open set $V\in\N_{bX}(x)$ such that $W=V\cap X$. There exists $a\in A\cap V$. We can find an ultrafilter $\U_a$ on $X$ such that $\lm_{bX} \U_a=a$. As $V$ is open, there exists $U\in\U_a$ such that $U\inc V\cap X=W$. Hence, 
$\bigwedge_{a\in A} \U_a\mesh\N_X(x)$. By (3), there 
exists 
$(a_n)_{n\in\omega}$ in $A$ such that 
$\bigwedge_{n\in \omega} \U_{a_n}\mesh\N_X(x)$. We claim 
that 
$x\in \cl_{bX}\{a_n:n\in\omega\}$. Indeed, for each $V\in \N_{bX}(x)$, we 
can find $V_1\in\N_{bX}(x)$ such that $\cl_{bX}V_1\inc V$. As $X\cap 
V_1\in\N_X(x)$ and $\N_X(x) \mesh \bigwedge_{n\in\omega} \U_{a_n}$, there 
is an $n$ such that $(X\cap V_1)\in\U_{a_n}$. Now, $a_n\in\cl_{bX}V_1\inc V$ and $x$ is an absolutely tight point.
\end{proof}
In view of Corollary \ref{cor:tightkernel}, we obtain:
\begin{cor}\cite{Ar.spectrum}
An absolutely tight point is productively countably tight.
\end{cor}
The converse is not true in general. Indeed, as observed in \cite{Ar.spectrum}, a $\Sigma$-product of uncountably many copies of the discrete two point space $\{0,1\}$ (or of the real line) is an example of a productively countably tight space (all points of which are, in particular, tight points) which is not absolutely countably tight.
However, the converse is true for countable spaces because every filter on a countable set is countably tight. More precisely:
\begin{prop}
Let $X$ be a topological space. The following are equivalent:
\begin{enumerate}
\item $X$ is productively countably tight;
\item $X$ is steadily countably tight and every countable subset of $X$ is productively countably tight;
\item $X$ is steadily countably tight and every countable subset of $X$ is absolutely tight.
\end{enumerate}
\end{prop}
\begin{proof}
The equivalence between (1) and (2) follows from (1) $\iff$ (2) in Corollary \ref{cor:tightkernel}.
Moreover, every filter on a countable set is countably tight, so that productive countable tightness and absolute tightness coincide on countable sets; and (2) $\iff$ (3) follows.
\end{proof}

\cite[Example 1]{bella.malykhin} is an example under (CH) of a tight point
  which does not have countable absolute tightness.  The space considered is countable, so that it provides an example of 
  a tight point which is not  productively countably
tight. The associated filter is in $\ker(\mathbb
{F}_1^{\bigtriangleup\diamondsuit})\setminus \ker(\mathbb
{F}_1^{\diamondsuit\diamondsuit})$. It shows that, at least under (CH), the converse of Corollary \ref{cor:prodtighttotight} is not true.

\subsection{Productively Fr\'echet and $\aleph_0$-bisequential spaces}

In \cite{Ar.spectrum} a regular space $X$ is called
$\aleph_0$-bisequential provided that $\omega$ is in the frequency
spectrum of $X$ (by \cite[Theorem 3.6]{Ar.spectrum} and
Corollary~\ref{cor:tight}, it means that $X$ is productively countably tight) and every countable subset of
$X$ is bisequential (\footnote{A filter $\F$ is {\em bisequential}
if for every filter $\G\mesh\F$ there exists a countably based filter
$\H\mesh\G$ such that $\H\geq \F$. A topological space is {\em
bisequential} if all its neighborhood filters are bisequential.}).
  \cite[Proposition 6.27]{Ar.spectrum} states that the product of an
 $\aleph_0$-bisequential space with a strongly Fr\'echet space is
 strongly Fr\'echet.   In view of Proposition \ref{th:mainJM}, every $\aleph_0$-bisequential space is productively Fr\'echet. We can also give a direct proof of this fact by characterizing productive Fr\'echetness in terms comparable to $\aleph_0$-bisequentiality.
\begin{thm}\label{thm:1}
$\F\in\mathbb {F}_{\omega}^{\bigtriangleup\bigtriangleup}(X)$ if
and only if $\F\in\mathbb
{F}_{\omega}^{\bigtriangleup\diamondsuit}(X)$ and $\F\vee
A\in\mathbb {F}_{\omega}^{\bigtriangleup\bigtriangleup}(X)$ for
every countable $A\subseteq X$ such that $A\mesh\F$.
\end{thm}

\begin{proof}
Suppose $\F\in\mathbb
{F}_{\omega}^{\bigtriangleup\diamondsuit}(X)$ and
 $\F\vee A\in\mathbb {F}_{\omega}^{\bigtriangleup\bigtriangleup}$ for every countable
 $A\subseteq X$ such that $A\mesh\F$.  Let $\G\in\mathbb {F}_{\omega}^{\bigtriangleup}$
 and $\G\mesh\F$.  Since $\F\in\mathbb {F}_{\omega}^{\bigtriangleup\diamondsuit}(X)$,
there is a countable set $A$ such that $A\mesh\G\vee\F$. Since
$\F\vee A\in\mathbb {F}_{\omega}^{\bigtriangleup\bigtriangleup}$,
and $\G\mesh(\F\vee A)$, there is countably based filter
$\C\geq\G\vee(\F\vee A) \geq\G\vee\F$.  Thus,
$\F\bigtriangleup\G$.  Therefore,
 $\F\in\mathbb {F}_{\omega}^{\bigtriangleup\bigtriangleup}$.

The opposite implication is trivial.
\end{proof}

We now deduce
\begin{cor}\cite[Proposition 6.27]{Ar.spectrum}\label{cor:1}
If $X$ is $\aleph_0$-bisequential, then $X$ is productively
Fr\'echet.
\end{cor}
\begin{proof}
Let $\F$ be neighborhood filter of a point in $X$. By
\cite[Theorem 3.6 ]{Ar.spectrum},
 $\F\in\ker(\mathbb {F}_{1}^{\diamondsuit\diamondsuit})\subseteq\mathbb {F}_{1}^{\diamondsuit\diamondsuit}$.
In particular, $\F\in\mathbb
{F}_1^{\diamondsuit\diamondsuit}\subseteq\mathbb
{F}_{\omega}^{\bigtriangleup\diamondsuit}$.
 Also, for every countable set $A$ meshing $\F$ we have $\F\vee A$ bisequential
 and hence in $\mathbb {F}_{\omega}^{\bigtriangleup\bigtriangleup}$.
So by Theorem~\ref{thm:1}, $\F\in \mathbb
{F}_{\omega}^{\bigtriangleup\bigtriangleup}$.
\end{proof}
There are models of set theory in which the converse of
Corollary~\ref{cor:1} is false \cite{JM}. It is unknown if there
is a ZFC example of a productively Fr\'echet space which is not
$\aleph_0$-bisequential.


\subsection{$\mathbb J$-expandable presentations and kernels of $\diamondsuit$-polars}
To apply Theorem \ref{thm:generaltightpolar} for a class $\mathbb J$ that is not stable under infima, we need to use general presentations of filters and not only infima as in Corollary \ref{cor:tightgeneral}. Practically speaking, to apply Theorem \ref{thm:generaltightpolar} in such cases, we need to characterize $\mathbb J$-expandability of a presentation.

A {\em crossing of a presentation} $\gamma:I\times J\to 2^X$ is a family $\{D_{\rho,l}:(\rho,l)\in I\times K\}$ such that 
 for every $\beta\in J$ and $l\in K$ there exists $\alpha\in I$ such that $\gamma(\alpha,\beta)\cap D_{\alpha,l}\neq\emptyset$. A crossing $\{D_{\rho,l}\colon (\rho,l)\in I\times\omega\}$ satisfying 
$D_{\rho,l+1}\subseteq D_{\rho,l}$ for every $(\rho,l)\in I\times\omega$ is called $\omega$-{\em decreasing}.
\begin{lem}\label{lem:pp1} 
A presentation $\gamma:I\times J\to 2^X$ is  $\Fo_{\omega}^{\bigtriangleup}$-expandable
if and only if for every $\omega$-decreasing crossing $\{D_{\rho,l}\colon (\rho,l)\in I\times\omega\}$ of $\gamma$, 
 there exist $e_k\in X$ and $\alpha_k\in I$, such that for every $\beta\in J$ and $n\in\omega$ we have $e_k\in\gamma(\alpha_k,\beta)\cap D_{\alpha_k,n}$ for all $k$ sufficiently large.
\end{lem}
\begin{proof}
Suppose $\gamma:I\times J\to X$ is not $\Fo_{\omega}^{\bigtriangleup}$-expandable.  
There is a decreasing countably based filter $\E=\{E_k\}_{k\in\omega}$  on $X\times I$ such that $\E\mesh\gamma^*$ and 
there is no sequence $(x_n,\alpha_n)\geq\E\vee\gamma^*$.  Since $\E\mesh\gamma^*$, for every $\beta\in J$ and $k\in\omega$ there is an $\alpha\in I$ such that
\begin{equation}\label{inproof1}
\gamma(\alpha,\beta)\times\{\alpha\}\cap E_k\neq\emptyset.
\end{equation}
For each $k\in\omega$ and $\alpha\in I$, let $D_{\alpha,k}=\{x\in X:(x,\alpha)\in E_k\}$.  By (\ref{inproof1}), $\{D_{\alpha,k}:(\alpha,k)\in I\times\omega\}$ is a crossing of $\gamma$. It is $\omega$-decreasing because $\{E_k:k\in\omega\}$ is decreasing.  By way of contradiction, assume that there exist points $e_l\in X$ and $\alpha_k\in I$, such that for every $\beta\in J$ and $k\in\omega$ we have $e_l\in\gamma(\alpha_l,\beta)\cap D_{\alpha_l,k}$ for all $l$ sufficiently large.  For each $l\in\omega$ let $y_l=(e_l,\alpha_l)$.  We show that $(y_l)_{l\in\omega}\geq\gamma^*\vee\E$.  Let $\beta\in J$ and $k\in\omega$.  There is an $n\in\omega$ such that $e_l\in\gamma(\alpha_k,\beta)\cap D_{\alpha_k,k}$ for all 
$l\geq n$.  Now $y_l=(e_l,\alpha_l)\in(\gamma(\alpha_l,\beta)\cap D_{\alpha_l,k})\times\{\alpha_l\}=(\gamma(\alpha_l,\beta)\times\{\alpha_l\})\cap E_k$ for all $l\geq n$.  So, $(y_l)_{l\in\omega}\geq\gamma^*\vee\E$, contradicting our choice of $\E$.

Suppose $\gamma:I\times J\to 2^X$ is an $\Fo_{\omega}^{\bigtriangleup}$-expandable presentation.  Let 
$\{D_{k,\alpha}\colon (k,\alpha)\in \omega\times I\}$ be an $\omega$-decreasing crossing of $\gamma$.  For each $k\in\omega$ define 
$E_k=\{(x,\alpha)\colon x\in D_{k,\alpha}\}$.  It is easily verified that $(E_k)_{k\in\omega}$ is a countably based filter meshing with $\gamma^*$.  So, there exists $(y_k)_{\omega}\geq(E_k)_{\omega}\vee\gamma^*$.  For each $k$ let $e_k=\pi_X(y_k)$ and $\alpha_k=\pi_{I}(y_k)$.  Let $\beta\in J$ and $n\in\omega$.  There is a $l\in\omega$ such that $y_k\in E_n\cap\gamma(\alpha_k,\beta)\times{\alpha_k}$ for all $k\geq l$.  So, $e_k\in D_{\alpha_k,n}\cap\gamma(\alpha_k,\beta)$ for all $k\geq l$.  
\end{proof}
\begin{lem}\label{lem:pp2}
A presentation $\gamma:I\times J\to 2^X$ is  $\Fo_{1}^{\dagger}$-expandable 
 if and only if for any crossing $\{D_{\rho,l}\colon (\rho,l)\in I\times\omega\}$ of $\gamma$, we may find finite sets $E_k\subseteq X$ and $H_k\subseteq I$, such that for every $\beta\in J$ there is  $k\in\omega$ and $\alpha\in H_k$ such that we have $E_k\cap\gamma(\alpha,\beta)\cap D_{\alpha,k}\neq\emptyset$. 
\end{lem}
\begin{proof}
Suppose $\gamma:I\times J\to 2^X$ is not $\Fo_{1}^{\dagger}$-expandable.  
There is a countable collection of sets $\{B_k\}_{k\in\omega}$ on $X\times I$ meshing with $\gamma^*$ such that
for any selection of finite sets $G_k\subseteq B_k$, $\bigcup_{k\in\omega}G_k$ does not mesh with 
$\gamma^*$.  Since $B_k\mesh\gamma^*$, for every $\beta\in J$ and $k\in\omega$ there is an $\alpha\in I$ such that
\begin{equation}\label{bun1}
\gamma(\alpha,\beta)\times\{\alpha\}\cap B_k\neq\emptyset.
\end{equation}
For each $k\in\omega$ and $\alpha\in I$ let $D_{\alpha,k}=\{x\in X: (x,\alpha)\in B_k\}$.  By (\ref{bun1}), it is a crossing of $\gamma$.  By way of contradiction, assume that there exist finite 
sets $E_k\in X$ and $H_k\in I$, such that for every $\beta\in J$ there is a $k\in\omega$ and a $\alpha\in H_k$ such that  $E_k\cap\gamma(\alpha,\beta)\cap D_{\alpha,k}\neq\emptyset$.  For each $k\in\omega$ let 
$G_k=\{(x,\alpha)\colon\alpha\in H_k\text{ and }x\in D_{\alpha,k}\cap E_k\}$.  Notice that $G_k$ is a 
finite subset of $B_k$.  Let $\beta\in J$.  There is a $k\in\omega$ and an $\alpha\in H_k$ such that 
$E_k\cap\gamma(\alpha,\beta)\cap D_{\alpha,k}\neq\emptyset$.  Now $G_k\cap(\gamma(\alpha,\beta)\times\{\alpha_l\})\neq\emptyset$.  So, $\bigcup_{k\in\omega}G_k\mesh\gamma^*$, 
contradicting our choice of $\{B_k\}_{k}$.

Suppose that $\gamma$ is $\Fo_{1}^{\dagger}$-expandable.  Let 
$\{D_{\alpha,k}\colon (\alpha,k)\in I\times\omega\}$ be a crossing of $\gamma$.  For each $k\in\omega$ define 
$B_k=\{(x,\alpha)\colon x\in D_{\alpha,k}\}$.  It is easily verified that $B_k\mesh\gamma^*$ for every 
$k\in\omega$.  So, there exist finite sets $G_k\subseteq B_k$ such that $\bigcup_{k\in\omega}G_k\mesh\gamma^*$.  For each $k$ let $E_k=\pi_X(G_k)$ and $H_k=\pi_{I}(G_k)$.  Let $\beta\in J$.  There is a $k\in\omega$ 
and an $\alpha\in I$ such 
that $G_k\cap(\gamma(\alpha,\beta)\times\{\alpha\})\neq\emptyset$.  Let $(x,\alpha)\in G_k\cap(\gamma(\alpha,\beta)\times\{\alpha\})$.  Notice that $\alpha\in H_k$ and $x\in E_k$.  Thus,  $E_k\cap D_{\alpha,k}\cap\gamma(\alpha,\beta)\neq\emptyset$ for some $\alpha\in H_k$.  
\end{proof}

Lemma \ref{lem:pp1} combined with Theorem \ref{thm:generaltightpolar} leads to:
\begin{thm}\label{thm:202}
The following are equivalent for a filter $\F\in\mathbb {F}(X)$.
\begin{enumerate}
\item $\F\in\ker(\mathbb {F}_{\omega}^{\bigtriangleup\diamondsuit})$;
\item $\F\in \mathbb {F}_{\omega}^{\bigtriangleup\diamondsuit}$ and $A\vee\F\in\ker(\mathbb {F}_{\omega}^{\bigtriangleup\diamondsuit})$ for all countable set $A$ meshing with $\F$;
\item for every presentation $\gamma:I\times J\to 2^X$ for which there exists points $e_k$ in $X$ and indices $\alpha_k\in I$ such that for every $\beta\in J$ and every $n\in \omega$ we have $e_k\in\gamma(\alpha,\beta)\cap D_{\alpha_k,n}$ for all $k$ sufficiently large, whenever 
$\{D_{\rho,l}\colon (\rho,l)\in I\times\omega\}$ is an $\omega$-decreasing crossing of $\gamma$: if $\gamma_*\mesh\F$, then there exists a countable subset $C$ of $I$ such that 
$\{\bigcup_{\alpha\in C} \gamma(\alpha,\beta)\colon\beta\in 
J\}\mesh\F$.
\end{enumerate}
\end{thm}

Theorem~\ref{thm:202} gives a concrete characterization of the second property in the following instance of Corollary \ref{cor:product}.
\begin{cor}
The following are equivalent for a topological space $X$.
\begin{enumerate}
\item the product of $X$ with every strongly Fr\'echet space is countably tight;
\item $X$ is  $\ker(\mathbb {F}_{\omega}^{\bigtriangleup\diamondsuit})$-based.
\end{enumerate}
\end{cor}

Lemma \ref{lem:pp2} combined with Theorem \ref{thm:generaltightpolar} leads to:
\begin{thm}\label{thm:203}
The following are equivalent for a filter $\F\in\mathbb {F}(X)$.
\begin{enumerate}
\item $\F\in\ker(\mathbb {F}_{1}^{\dagger\diamondsuit})$.
\item $\F\in \mathbb {F}_{1}^{\dagger\diamondsuit}$ and $A\vee\F\in\ker(\mathbb {F}_{1}^{\dagger\diamondsuit})$ for all countable set $A$ meshing with $\F$;
\item for any presentation $\gamma:I\times J\to 2^X$ for which there exist finite sets $E_k\subseteq X$ and $H_k\subseteq I$, such that for every $\beta\in J$ there is  $k\in\omega$ and $\alpha\in H_k$ verifying $E_k\cap\gamma(\alpha,\beta)\cap D_{\alpha,k}\neq\emptyset$ whenever $\{D_{\rho,l}\colon (\rho,l)\in I\times\omega\}$ is a crossing of $\gamma$: if $\gamma_*\mesh\F$, then there is a countable subset $C$ of $I$ such that 
$\{\bigcup_{\alpha\in C} \gamma(\alpha,\beta)\colon\beta\in 
J\}\mesh\F$.
\end{enumerate}
\end{thm}
Theorem~\ref{thm:203} gives a concrete characterization of the second property in the following instance of Corollary \ref{cor:product}.
\begin{cor}
The following are equivalent for a topological space $X$.
\begin{enumerate}
\item the product of $X$ with every countably fan-tight space is countably tight;
\item $X$ is  $\ker(\mathbb {F}_{1}^{\dagger\diamondsuit})$-based.
\end{enumerate}
\end{cor}
Characterizing kernels of $\dagger$-polars internally seems to require even more machinery. While we do not have concrete characterizations of these kernels, we can give some information (in the next section) 
on how they relate to kernels of $\diamondsuit$-polars that we have characterized before.

\section{Inclusions}

\begin{thm}\label{thm:5}
 $\ker(\mathbb
{F}_{1}^{\bigtriangleup\diamondsuit})\subseteq\ker(\mathbb
{F}_{\omega}^{\bigtriangleup\dagger})$
\end{thm}
\begin{proof}

Let $\F\in\ker(\mathbb {F}_{1}^{\bigtriangleup\diamondsuit})$ on
$X$, $Y$ be a set, and $A\subseteq X\times Y$.  We show that
$A\F\in\mathbb {F}_{\omega}^{\bigtriangleup\dagger}(Y)$, using Lemma \ref{lem:fanT}.  Let
$\G\in\mathbb {F}_{\omega}^{\bigtriangleup}(Y)$ and
$\B=\{B_k\}_{k\in\omega}$ be a countably based filter such that
$B_k\mesh(\G\vee A\F)$ for every $k\in\omega$.

Fix $F\in\F$.  Since $\G\in\mathbb
{F}_{\omega}^{\bigtriangleup}(Y)$ and $\G\mesh(AF\vee\B)$, we may
find a sequence $(y^{F}_{k})_{k\in\omega}$ such that
$(y^{F}_{k})_{k\in\omega}\geq AF\vee\G$ and $y^F_k\in B_k$ for
every $k\in\omega$.  For each $k\in\omega$ define $x_{k}^F\in F$
so that $(x^F_k,y^F_k)\in A$.

Since $\F\in\ker(\mathbb {F}_{1}^{\bigtriangleup\diamondsuit})$ on
$X$, and $\F\mesh\bigwedge_{F\in\F}(x^{F}_k)_{k\in\omega}$, we may
find $\{F_n\colon\in\omega\}\subseteq\F$ such that
$\F\mesh\bigwedge_{n\in\omega}(x^{F_n}_k)_{k\in\omega}$.

For every $k\in\omega$ let $T_k=\{y^{F_n}_k\colon n\leq k\}$.
Notice that $T_k\subseteq B_k$ and $T_k$ is finite for all
$k\in\omega$.  Let $F\in\F$ and $G\in\G$.  There is an
$n\in\omega$ such that $(x^{F_n}_k)_{k\in\omega}\mesh F$.  Since
$(y^{F_n}_k)_{k\in\omega}\geq\G$ and $x^{F_n}_k\in F$ for
infinitely many $k\in\omega$, there is a $p\geq n$ such that
$x^{F_n}_p\in F$ and $y^{F_n}_p\in G$.  Hence, $y^{F_n}_p\in
AF\cap G$.  Since $n\leq p$, we have $T_p\cap AF\cap
G\neq\emptyset$.  So, $\bigcup_{k\in\omega}T_k\mesh(\G\vee A\F)$.
Thus, $A\F\in\mathbb {F}_{\omega}^{\bigtriangleup\dagger}(Y)$.
\end{proof}
This improves significantly  \cite[Proposition 2.1]{bella.vanmill} that states the weaker inclusion $\ker(\mathbb F_1^{\bigtriangleup\diamondsuit})\inc \mathbb F_1^\dagger$  in topological terms (i.e., every tight point has countable fan-tightness) under the assumption of $T_1$.

\begin{cor}
The product of a space whose every point is tight with a strongly Fr\'echet space has countable fan-tightness.
\end{cor}
\begin{thm}\label{thm:6}
 $\ker(\mathbb
{F}_{1}^{\diamondsuit\diamondsuit})\subseteq\ker(\mathbb
{F}_{1}^{\dagger\dagger})$.
\end{thm}
\begin{proof}

Let $\F\in\ker(\mathbb {F}_{1}^{\diamondsuit\diamondsuit})$ on
$X$, $Y$ be a set, and $A\subseteq X\times Y$.  We show that
$A\F\in\mathbb {F}_{1}^{\dagger\dagger}(Y)$.  Let $\G\in\mathbb
{F}_{1}^{\dagger}(Y)$ and $\B=\{B_k\}_{k\in\omega}$ be a countably
based filter such that $B_k\mesh(\G\vee A\F)$ for every $k\in\omega$.

Fix $F\in\F$.  Since $\G\in\mathbb {F}_{1}^{\dagger}(Y)$ and
$\G\mesh(AF\vee\B)$, we may find finite sets $C^F_k\subseteq B_k\cap
AF$ such that $(\bigcup_{k\geq n}C_k)\mesh AF\vee\G$ for every
$n\in\omega$. Let $\H_{F}=(\bigcup_{k\geq
n}C^F_k)_{n\in\omega}\vee\G$.  In view of Corollary \ref{cor:supfantight}, $\H_F\in\mathbb
{F}_1^{\dagger}(Y)$.  Define a function
$f^F\colon\bigcup_{k\in\omega}C^F_k\to F$ such that $f^F\subseteq A$
and define $\G_{F}=f^F(\H_F)$.

Since $\F\in\ker(\mathbb {F}_{1}^{\diamondsuit\diamondsuit})$ on
$X$, $\G_F\in\mathbb {F}_1^{\diamondsuit}$ for all $F\in\F$, and
$\F\mesh\bigwedge_{F\in\F}\G_F$, we may find
$\{F_n\colon\in\omega\}\subseteq\F$ such that
$\F\mesh\bigwedge_{n\in\omega}\G_{F_n}$.

For every $k\in\omega$ let $T_k=\bigcup\{C^{F_n}_k\colon n\leq
k\}$.  Notice that $T_k\subseteq B_k$ and $T_k$ is finite for all
$k\in\omega$.  Let $F\in\F$ and $G\in\G$.  There is an
$n\in\omega$ such that $\G_{F_n}\mesh F$.  Thus,
$f^{F_n}\left[G\cap\bigcup_{k\geq n}C^{F_n}_k\right]\cap
F\neq\emptyset$.  So, there is a $k\geq n$ such that
$f^{F_n}[G\cap C^{F_n}_k]\cap F\neq\emptyset$.  Since
$f^{F_n}\subseteq A$, we have $G\cap C^F_k\cap AF\neq\emptyset$.
So, $\bigcup_{k\in\omega}T_k\mesh(\G\vee A\F)$.  Thus, $A\F\in\mathbb
{F}_{1}^{\dagger\dagger}(Y)$.
\end{proof}
\cite[Theorem 5]{ar.bella.CMUC} proves the much weaker inclusion $\ker(\mathbb F_1^{\diamondsuit\diamondsuit})\inc \mathbb F_1^\dagger$ in topological terms (i.e., if the product of $X$ with every space of countable tightness has countable tightness, then $X$ has countable fan-tightness) among Tychonoff spaces. From the theorem above, we conclude:
\begin{cor}\label{cor:tofantight}
If $X$ is productively countably tight (in particular if it has countable absolute tightness), then its product with every space of countable fan-tightness has countable fan-tightness.
\end{cor}
This result generalizes \cite[Corollary 6]{ar.bella.CMUC} that states the same conclusion (among regular spaces) under the assumption that $X$ is a compact space of countable tightness, which, of course, implies that $X$ has countable absolute tightness.

By definition, $\mathbb {F}_1^{\diamondsuit\diamondsuit}\subseteq
\mathbb {F}_1^{\bigtriangleup\diamondsuit}$, so that $\ker(\mathbb
{F}_1^{\diamondsuit\diamondsuit})\subseteq \ker(\mathbb
{F}_1^{\bigtriangleup\diamondsuit})$.  By Theorem~\ref{thm:5},
$\ker(\mathbb {F}_1^{\bigtriangleup\diamondsuit})\subseteq \mathbb
{F}_1^{\dagger}$.  Therefore, $\ker(\mathbb
{F}_1^{\diamondsuit\diamondsuit})\subseteq \ker(\mathbb
{F}_1^{\bigtriangleup\diamondsuit})\subseteq \mathbb
{F}_1^\dagger\subseteq\mathbb {F}_1^\diamondsuit$.
All these inclusions are strict in general. For instance,
$\S_\omega\in \mathbb {F}_1^\diamondsuit\setminus\mathbb
{F}_1^\dagger$. Moreover, we already have observed that
 \cite[Example 1]{bella.malykhin} shows that  $\ker(\mathbb
{F}_1^{\bigtriangleup\diamondsuit})\setminus \ker(\mathbb
{F}_1^{\diamondsuit\diamondsuit})\neq\varnothing$ under (CH). Finally, \cite[Example
3]{bella.malykhin} is a ZFC example of a point of countable
fan-tightness which is not tight. The neighborhood filter is in
$\mathbb {F}_1^\dagger \setminus \ker(\mathbb
{F}_1^{\bigtriangleup\diamondsuit})$.  Generally, we have the
following picture showing the containments that we know.

\bigskip

\begin{picture}(300,130)

\thicklines \put(10,50){\makebox(0,0){$\ker(\mathbb
{F}_1^{\diamondsuit\diamondsuit})$}}
\put(10,90){\makebox(0,0){$\ker(\mathbb
{F}_{1}^{\bigtriangleup\diamondsuit})$}}
\put(150,90){\makebox(0,0){$\ker(\mathbb
{F}_{\omega}^{\bigtriangleup\dagger})$}}
\put(150,50){\makebox(0,0){$\ker(\mathbb
{F}_{\omega}^{\bigtriangleup\diamondsuit})$} }
\put(10,10){\makebox(0,0){$\ker(\mathbb
{F}_{1}^{\dagger\dagger})$}}
\put(150,10){\makebox(0,0){$\ker(\mathbb
{F}_{1}^{\dagger\diamondsuit})$} }
\put(240,90){\makebox(0,0){$\mathbb {F}_{1}^{\dagger}$}}
\put(240,50){\makebox(0,0){$\mathbb
{F}_{1}^{\diamondsuit}$}} 
\put(150,130){\makebox(0,0){$\mathbb
{F}_{\omega}^{\bigtriangleup\bigtriangleup}$}}
\put(240,130){\makebox(0,0){$\mathbb
{F}_{\omega}^{\bigtriangleup}$}}

\put(240,80){\vector(0,-1){20}}
\put(175,90){\vector(1,0){40}}
\put(175,130){\vector(1,0){40}}
\put(10,60){\vector(0,1){20}}
\put(10,40){\vector(0,-1){20}}
\put(35,15){\vector(3,2){90}}
\put(175,50){\vector(1,0){40}}
\put(40,10){\vector(1,0){85}}
\put(40,90){\vector(1,0){85}}
\put(150,80){\vector(0,-1){20}}
\put(150,15){\vector(0,1){20}}
\put(150,120){\vector(0,-1){20}}
\put(240,120){\vector(0,-1){20}}

\end{picture}

\bigskip


\cite[Theorem 2.3]{bella.vanmill} states that in a regular countably compact space,
points of countable tightness are tight.
This result can be stated at the level of filters via
\eqref{eq:collapse}.
 If $X$ is a topological
space, we denote by $\O_X$ the class of filters admitting a base
composed of open sets and by $\K_\omega$ the class of filters admitting a base
of countably compact sets. The result quoted above follows from
\begin{equation}\label{eq:collapse}
\K_\omega\cap\O\cap\ker(\mathbb
{F}_1^{\bigtriangleup\diamondsuit})=
\K_\omega\cap\O\cap \mathbb {F}_1^\diamondsuit,
\end{equation}
by considering neighborhood filters.

 We show that the coincidence of these classes
 is true not only for filters of $\K_\omega$,
but for the broader class of strongly $q$-regular filters.
A filter $\F$ is {\em strongly} $q$-{\em regular} if
for every $F\in \F$, there exists a sequence $(Q_n^F)_{n\in\omega}$ in $\F$
such that
$\adh_X (\bigwedge_{i\in\omega} (x_n^i)_{n\in\omega})\cap 
F\neq\varnothing$
whenever $\bigwedge_{i\in\omega} (x_n^i)_{n\in\omega}\mesh(Q_n^F)_{n\in\omega}$.
Let $\Q_\omega$ denote the class of strongly $q$-regular filters.
We say that a point is strongly $q$-regular if its
 neighborhood filter is strongly $q$-regular.

This notion is a little bit stronger than that
of a $q$-regular point in the sense of \cite{dolnog}, which is a variant
of a regular $q$-point in the sense of \cite{quest}.
All points of countable character and all points of
a regular locally countably compact space
are strongly $q$-regular.
\cite[Proposition 13]{dolnog} gives an example of a non first-countable and non locally countably compact
strongly $q$-regular space.

\begin{thm}\label{th:coincide}
$$\Q_\omega\cap\O\cap\ker(\mathbb {F}_1^{\bigtriangleup\diamondsuit})=\Q_\omega\cap\O\cap \mathbb
{F}_1^\diamondsuit.$$
\end{thm}
\begin{proof}
Let $\F\in \Q_{\omega}\cap\O\cap\mathbb {F}_1^\diamondsuit$.
For every $F\in \F$ there exists $(Q_n^F)_{n\in\omega}$ witnessing
the definition of a strongly $q$-regular filter. Let
$\bigwedge_{\alpha\in I} (x_k^\alpha)_{k\in\omega}\mesh\F$. In
particular, $\bigwedge_{\alpha\in I} 
(x_k^\alpha)_{k\in\omega}\mesh Q_n^F$ for every $n\in \omega$ and
$F\in\F$. Thus, there exists
$(x_k^{\alpha_{n,F}})_{k\in\omega}\mesh Q_n^F$. Hence
$\bigwedge_{n\in\omega} (x_k^{\alpha_{n,F}})_{k\in\omega}
\mesh(Q_n^F)_{n\in\omega}$. Therefore, there exists $x_F\in F\cap
\adh \bigwedge_{n\in\omega} (x_k^{\alpha_{n,F}})_{k\in\omega}$.
Now $\{x_F: F\in \F\}\mesh\F$ and $\F\in\mathbb F_1^\diamondsuit$, so that there exists a sequence
$(F_i)_{i\in\omega}$ in $\F$ such that $\{x_{F_i}:
i\in\omega\}\mesh\F$. Then 
$\bigwedge_{n\in\omega, i\in\omega}
(x_k^{\alpha_{n,F_i}})_{k\in\omega}\mesh\F$. Indeed, for every open
$U\in \F$, there exists $i\in \omega$ such that $x_{F_i}\in U$.
But $x_{F_i}\in
 \adh \bigwedge_{n\in\omega}
(x_k^{\alpha_{n,F_i}})_{k\in\omega}$. As $U$ is open,
$\bigwedge_{n\in\omega}
(x_k^{\alpha_{n,F_i}})_{k\in\omega}\mesh U$.
\end{proof}

The following generalizes \cite[Theorem 2.3]{bella.vanmill}:
\begin{cor} \label{cor:coincide}
A strongly $q$-regular point of countable tightness is tight.
\end{cor}


A filter $\F$ is called {\em regularly of pointwise countable type}
if for every $F\in \F$, there exists a sequence
$(Q_n^F)_{n\in\omega}$ in $\F$ such that
$\adh_X \H\cap F\neq\varnothing$ whenever
$\H\mesh(Q_n^F)_{n\in\omega}$. Let $\Q$ denote the class of regularly of pointwise countable
type filters.

Adapting the proof of Theorem \ref{th:coincide} and letting $\mathbb{A}$ stand for the class of 
absolutely countably tight fiters, we get:

\begin{thm}
$$\Q\cap\O\cap\mathbb{A}=\Q\cap\O\cap\ker(\mathbb {F}_1^{\diamondsuit\diamondsuit})=\Q\cap\O\cap
\mathbb{F}_1^\diamondsuit.$$
\end{thm}
Notice that a regular point of pointwise countable type in the sense of \cite{quest}
has a neighborhood filter in $\Q$.
\begin{cor}\label{cor:ptwise}
A point  of pointwise countable type and of countable tightness of a completely regular space is
absolutely countably tight.
\end{cor}
Combined with Corollary \ref{cor:tofantight}, we obtain another improvement of \cite[Corollary 6]{ar.bella.CMUC}:
\begin{cor}
If $X$ is regular, of pointwise countable type and of countable tightness, then its product with every space of countable fan-tightness has countable fan-tightness.
\end{cor}
Corollary \ref{cor:ptwise} has other interesting consequences. For instance, it can be combined with Arhangel'skii 's theorem \cite[Theorem 1]{AVAM} stating that a topological group with an everywhere dense subset of countable absolute tightness is metrizable to the effect that
\begin{thm}
If a (completely regular) topological group has an everywhere dense subset that has countable tightness and is of pointwise countable type, then this group is metrizable.
\end{thm}

\bibliographystyle{plain}

\end{document}